\documentclass{amsart}

\usepackage{amssymb}
\usepackage{amsmath}
\usepackage{amsthm}
\usepackage{graphicx}

\theoremstyle{plain}
\newtheorem{theorem}{Theorem}[section]

\newtheorem{lemma}[theorem]{Lemma}
\newtheorem{corollary}[theorem]{Corollary}
\newtheorem{proposition}[theorem]{Proposition}

\theoremstyle{definition}
\newtheorem{definition}[theorem]{Definition}
\newtheorem{example}[theorem]{Example}

\theoremstyle{remark}
\newtheorem{remark}[theorem]{Remark}
\newtheorem*{notation}{Notation}
\newtheorem*{acknowledgments}{Acknowledgments}

\numberwithin{figure}{section}

\begin{document}

\title[Heegaard surfaces and measured laminations II]{Heegaard surfaces and measured laminations, II: non-Haken 3--manifolds}

\author{Tao Li}
\thanks{Partially supported by an NSF grant}

\address{Department of Mathematics \\
  Boston College \\
Chestnut Hill, MA 02467 \\
 USA}
\email{taoli@bc.edu}
\urladdr{http://www2.bc.edu/\~{}taoli}

\begin{abstract}
A famous example of Casson and Gordon shows that a Haken 3--manifold can have an infinite family of irreducible Heegaard splittings with different genera.
In this paper, we prove that a closed non-Haken 3--manifold has only finitely many irreducible Heegaard splittings, up to isotopy.  This is much stronger than the generalized Waldhausen conjecture.  Another immediate corollary is that for any irreducible non-Haken 3--manifold $M$, there is a number $N$, such that any two Heegaard splittings of $M$ are equivalent after at most $N$ stabilizations.
\end{abstract}

\maketitle

\setcounter{tocdepth}{1}
\tableofcontents

\section{Introduction}
A Heegaard splitting of a closed orientable 3--manifold is said to be reducible if there is an essential simple closed curve in the Heegaard surface bounding disks in both handlebodies.  Haken proved that a Heegaard splitting of a reducible 3--manifold is always reducible \cite{H}.  

The classification of irreducible Heegaard splittings has been a long-standing fundamental problem in 3--manifold topology.  Such classification has been achieved for certain non-hyperbolic manifolds, such as $S^3$ by Waldhausen \cite{W1}, Lens spaces by Bonahon and Otal \cite{BO}, and Seifert fiber spaces by \cite{BCZ, M, MSch}.  The main theorem of this paper is a finiteness result for non-Haken 3--manifolds.

\begin{theorem}\label{main}
A closed orientable non-Haken 3--manifold has only finitely many irreducible Heegaard splittings, up to isotopy.
\end{theorem}

An important question in the study of Heegaard splittings is whether there are ways to construct different Heegaard splittings.  By adding trivial handles, one can always construct an infinite family of Heegaard splittings for every 3--manifold.  Theorem~\ref{main} says that, for irreducible non-Haken manifolds, adding trivial handles is virtually the only way of obtaining new Heegaard splittings.

The study of Heegaard splitting has been dramatically changed since Casson and Gordon introduced the notion of strongly irreducible Heegaard splitting \cite{CG}.  They showed that \cite{CG} an irreducible Heegaard splitting of a non-Haken 3--manifold is also strongly irreducible.  Using the thin-position argument, Rubinstein established relations between strongly irreducible Heegaard splittings and normal surface theory.  The results in \cite{CG} have also been used to attack the virtually Haken conjecture \cite{La, MMZ}.  

Casson and Gordon found the first 3--manifolds containing infinitely many different irreducible Heegaard splittings, see \cite{CG2,Sed,Ko}, and Theorem~\ref{main} says that this can only happen in Haken 3--manifolds.  In section~\ref{Sexample}, we will show the relation between an incompressible surface and the infinite family of strongly irreducible Heegaard splittings in the Casson-Gordon example.  This interpretation of the Casson-Gordon example was independently discovered by \cite{MSS}, where the authors proved a special case of the theorem.

A conjecture of Waldhausen \cite{W2} says that a closed orientable 3--manifold has only finitely many minimal/reducible Heegaard splittings, up to homeomorphism (or even isotopy).  This is known to be false because of the Casson-Gordon example.  A modified version of this conjecture is the so-called generalized Waldhausen conjecture, which says that an irreducible and atoroidal 3--manifold has only finitely many Heegaard splittings in each genus, up to isotopy.  Johannson \cite{Jo1, Jo2} proved the generalized Waldhausen conjecture for Haken 3--manifolds.  Together with Johannson's theorem, Theorem~\ref{main} implies the generalized Waldhausen conjecture.  Moreover, Theorem~\ref{main} says that the original version of Waldhausen conjecture is true for non-Haken 3--manifolds.

Another important question in the study of Heegaard splittings is how different Heegaard splittings are related.  This is the so-called stabilization problem, asking the number of stabilizations required to make two Heegaard splittings equivalent.  It has been shown that the number of stabilizations is bounded by a linear function of the genera of the two splittings \cite{RS}, but it remains unknown whether there is a universal bound.  We hope the techniques used in this paper can shed some light on this question.  Corollary~\ref{stable} follows from Theorem~\ref{main} and \cite{RS}.  

\begin{corollary}\label{stable}
For any closed, orientable, irreducible and non-Haken 3--manifold $M$, there is a number $N$ such that any two Heegaard splittings of $M$ are equivalent after at most $N$ stabilizations.
\end{corollary}

We briefly describe the main ideas of the proof.  The basic idea is similar in spirit to the proof of \cite{L5}.  By \cite{H, CG, BCZ, BO, M, MSch}, we may assume $M$ is irreducible, atoroidal and not a small Seifert fiber space, and the Heegaard splittings are strongly irreducible.  By a theorem in \cite{L4}, there is a finite collection of branched surfaces in $M$ such that every strongly irreducible Heegaard surface is fully carried by a branched surface in this collection. Moreover, the branched surfaces in this collection have some remarkable properties, such as they do not carry any normal 2--sphere or normal torus.  Each surface carried by a branched surface corresponds to an integer solution to the system of branch equations \cite{FO}.  One can also define the projective lamination space for a branched surfaces, see \cite{L4}. If a branched surface in this collection carries an infinite number of strongly irreducible Heegaard surfaces, then we have an infinite sequence of points in the projective lamination space.  By compactness, there must be an accumulation point which corresponds to a measured lamination $\mu$.  The main task is to prove that $\mu$ is incompressible and hence yields a closed incompressible surface, contradicting the hypothesis that $M$ is non-Haken.  The proof utilizes properties of both strongly irreducible Heegaard splittings and measured laminations.

We organize this paper as follows.  In section~\ref{Sbranch}, we briefly review some results from \cite{L4} and show some relations between branched surfaces and strongly irreducible Heegaard splittings.  In sections~\ref{Slam} and \ref{Slimit}, we prove some technical lemmas concerning measured laminations.  In section~\ref{Shelix}, we explain a key construction.  We finish the proof of Theorem~\ref{main} in section~\ref{Smain}.  In section~\ref{Sexample}, we show how to interpret the limit of the infinite family of strongly irreducible Heegaard surfaces in the Casson-Gordon example. 

\begin{acknowledgments}
I would like to thank Bus Jaco, Saul Schleimer and Ian Agol for useful conversations and Cynthia Chen for technical assistance.  I also thank the referee for many corrections and suggestions.
\end{acknowledgments}

\section{Heegaard surfaces and branched surfaces}\label{Sbranch}

\begin{notation}
Throughout this paper, we will denote the interior of $X$ by $int(X)$, the closure (under path metric) of $X$ by $\overline{X}$, and the number of components of $X$ by $|X|$.  We will also use $|n|$ to denote the absolute value of $n$ if $n$ is a number.  We will use $\eta(X)$ to denote the closure of a regular neighborhood of $X$.  We will also use the same notations on branched surfaces and laminations as in sections 2 and 3 of \cite{L4}.
\end{notation}

Let $M$ be a closed orientable and non-Haken 3--manifold.  A theorem of Haken \cite{H} says that a reducible 3--manifold cannot have any irreducible Heegaard splitting.  By \cite{BCZ, BO, M, MSch}, Theorem~\ref{main} is true for small Seifert fiber spaces. So we may assume $M$ is irreducible and not a small Seifert fiber space.  Casson and Gordon \cite{CG} showed that irreducible Heegaard splittings are equivalent to strongly irreducible Heegaard splittings for non-Haken 3--manifolds.  Hence we assume the Heegaard splittings in this paper are strongly irreducible.  We call the Heegaard surface of a strongly irreducible splitting a strongly irreducible Heegaard surface.

By \cite{R, St}, each strongly irreducible Heegaard surface is isotopic to an almost normal surface with respect to a triangulation.  Similar to \cite{FO}, we can use normal disks and almost normal pieces to construct a finite collection of branched surfaces such that each strongly irreducible Heegaard surface is fully carried by a branched surface in this collection.  By a theorem of \cite{L4} (Theorem~\ref{Heeg1} below), we can split these branched surfaces into a larger collection of branched surfaces so that each strongly irreducible Heegaard surface is still fully carried by a branched surface in this collection and no branched surface in this collection carries any normal 2--sphere or normal torus. 

\begin{theorem}[Theorem 1.3 in \cite{L4}]\label{Heeg1}
Let $M$ be a closed orientable irreducible and atoroidal 3--manifold, and suppose $M$ is not a Seifert fiber space.  Then $M$ has a finite collection of branched surfaces, such that
\begin{enumerate}
\item each branched surface in this collection is obtained by gluing together normal disks and at most one almost normal piece with respect to a fixed triangulation, similar to \cite{FO},
\item up to isotopy, each strongly irreducible Heegaard surface is fully carried by a branched surface in this collection,
\item no branched surface in this collection carries any normal 2--sphere or normal torus. 
\end{enumerate}
\end{theorem}

Our goal is to prove that each branched surface in Theorem~\ref{Heeg1} only carries a finite number of strongly irreducible Heegaard surfaces.  We will use various properties of strongly irreducible Heegaard splittings, branched surfaces and measured laminations, and we refer to sections 2 and 3 of \cite{L4} for an overview of some results and techniques in these areas.  In this section, we prove some easy lemmas which establish some connections between branched surfaces and Heegaard surfaces.

Let $B$ be a branched surface, $N(B)$ be a fibered neighborhood of $B$, and $\pi:N(B)\to B$ be the map collapsing each $I$--fiber of $N(B)$ to a point. 
We say an annulus $A=S^1\times I\subset N(B)$ is a \emph{vertical annulus} if every $\{x\}\times I\subset A$ ($x\in S^1$) is a subarc of an $I$--fiber of $N(B)$.  We say a surface $\Gamma$ is carried by $N(B)$ if $\Gamma\subset N(B)$ is transverse to the $I$--fibers of $N(B)$.

\begin{figure}
\begin{center}
\includegraphics{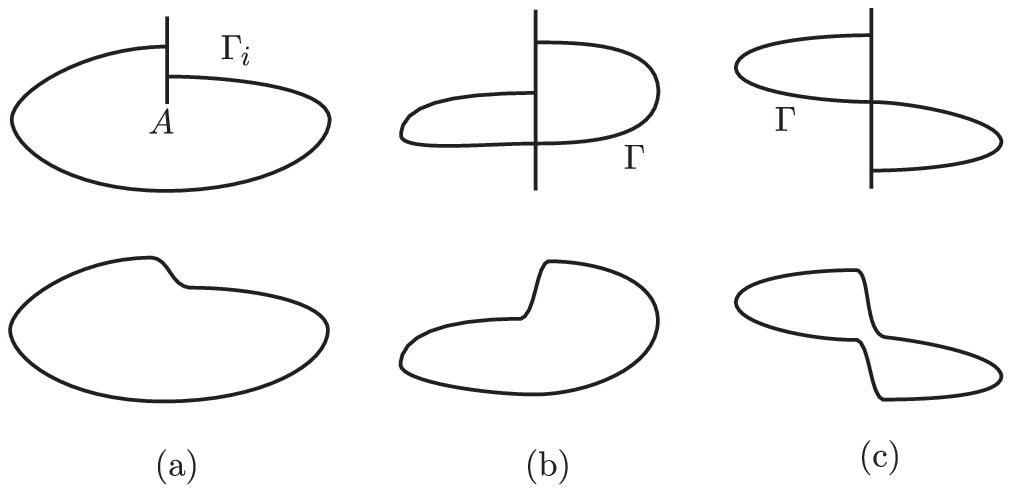}
\caption{}\label{tori}
\end{center}
\end{figure}

\begin{proposition}\label{Ptorus}
Let $B$ be a branched surface and $A\subset N(B)$ an embedded vertical annulus.  Suppose there is an embedded annulus $\Gamma$ carried by $N(B)$ such that $\partial\Gamma\subset A$ and $int(\Gamma)\cap A$ is an essential closed curve in $\Gamma$.  Then $B$ carries a torus.
\end{proposition}
\begin{proof}
First note that if $B$ carries a Klein bottle $K$, then the boundary of a twisted $I$--bundle over $K$ is a torus carried by $B$.  The idea of the proof is that one can perform some cutting and pasting on $A$ and $\Gamma$ to get a torus (or Klein bottle) carried by $B$.  The circle $int(\Gamma)\cap A$ cuts $\Gamma$ into 2 sub-annuli, say $\Gamma_1$ and $\Gamma_2$, with $int(\Gamma_i)\cap A=\emptyset$ ($i=1,2$).  Let $A_i$ be the sub-annulus of $A$ bounded by $\partial\Gamma_i$.  So $A_i\cup\Gamma_i$ is an embedded torus (or Klein bottle).  We have two cases here.  The first case is that $\Gamma_i$ connects $A$ from different sides, more precisely, after a small perturbation, the torus (or Klein bottle) $A_i\cup\Gamma_i$ is transverse to the $I$--fibers of $N(B)$, as shown in Figure~\ref{tori} (a).  The second case is that both $\Gamma_1$ and $\Gamma_2$ connect $A$ from the same side.  Then as shown in Figure~\ref{tori}(b, c), we can always use the annuli $\Gamma_i$ and $A_i$ to assemble a torus (or Klein bottle) carried by $B$.
\end{proof}

The following lemma is a variation of Lemma 2.2 in \cite{S} and the proof is  similar.

\begin{lemma}\label{Lnest}
Let $M=H_1\cup_SH_2$ be a strongly irreducible Heegaard splitting, $S$  the Heegaard surface, and $D$  an embedded disk in $M$ with $\partial D\subset S$.  Suppose $D$ is transverse to $S$ and $int(D)\cap S$ is a single circle $\gamma$.  Let $D_1\subset D$ be the disk bounded by $\gamma$, and suppose $D_1\subset H_1$ is a compressing disk of the handlebody $H_1$.  Then the annulus $A=D-int(D_1)$ must be $\partial$--parallel in the handlebody $H_2$.
\end{lemma}
\begin{proof}
Since $S$ is strongly irreducible, $\gamma$ does not bound a disk in $H_2$. So $A$ is incompressible in $H_2$, and hence $A$ is $\partial$--compressible.  Let $E\subset H_2$ be a $\partial$--compressing disk for the annulus $A=D-int(D_1)$.  We may suppose $\partial E$ consists of two arcs, $\alpha$ and $\beta$, where $\alpha\subset A$ is an essential arc in $A$, $\beta\subset S$ and $\partial\alpha=\partial\beta\subset\partial A$.  

Now we compress $A$ along $E$, in other words, we perform a simple surgery, replacing a small neighborhood of $\alpha$ in $A$ by two parallel copies of $E$.  The resulting surface is a disk properly embedded in $H_2$. We denote this disk by $D_2$.  After a small perturbation, we may assume $\partial D_2$ is disjoint from $\partial D_1$. Since $M=H_1\cup_SH_2$ is a strongly irreducible Heegaard splitting and $D_1$ is a compressing disk in $H_1$, $D_2$ must be a $\partial$--parallel disk in $H_2$.  This implies that $A$ is $\partial$--parallel in $H_2$.
\end{proof}

The following lemma follows easily from Proposition~\ref{Ptorus} and Lemma~\ref{Lnest}.

\begin{lemma}\label{Lbound}
Let $S$ be a strongly irreducible Heegaard surface fully carried by a branched surface $B$, and suppose $B$ does not carry any torus.  Let $A$ be an embedded vertical annulus in $N(B)$, and suppose $A\cap S=\cup_{i=1}^nc_i$ consists of $n$ non-trivial circles in $S$.  If some $c_i$ bounds a compressing disk in one of the two handlebodies, then there is a number $K$ depending only on $B$ such that $n=|A\cap S|<K$.
\end{lemma}
\begin{proof}
Suppose $M=H_1\cup_SH_2$ is the Heegaard splitting.  Let $A_i$ be the sub-annulus of $A$ bounded by $c_i\cup c_{i+1}$, and we may assume $A_i$ is properly embedded in $H_1$ if $i$ is odd and in $H_2$ if $i$ is even.  Without loss of generality, we may suppose $c_1$ bounds a compressing disk in a handlebody.  Note that the argument works fine if one starts with an arbitrary $c_i$ rather than $c_1$.

If $c_1$ bounds a compressing disk in $H_2$, since $c_1\cup c_2$ bounds an annulus $A_1$ in $H_1$, by Lemma~\ref{Lnest}, $A_1$ is $\partial$--parallel in $H_1$.  By pushing $A_1$ into $H_2$, we have that $c_2$ bounds a disk in $H_2$.  Since $A_2$ lies in $H_2$, the union of $A_2$ and the disk bounded by $c_2$ in $H_2$ is a disk bounded by $c_3$.  Since each $c_i$ is non-trivial in $S$, $c_3$ bounds a compressing disk in $H_2$.   Again, since $A_3$ lies in $H_1$, by Lemma~\ref{Lnest}, $A_3$ is $\partial$--parallel in $H_1$.  Inductively, we conclude that $A_{2k+1}$ is $\partial$--parallel in $H_1$ for each $k$.  So for each $k$, there is an annulus $\Gamma_k\subset S$ such that $\partial\Gamma_k=\partial A_{2k+1}$ and $A_{2k+1}\cup\Gamma_k$ bounds a solid torus $T_k$ in $H_1$.  It is clear that any two such solid tori $T_i$ and $T_j$ are either disjoint or nested.

Suppose $T_i$ and $T_j$ are nested, say $T_i\subset T_j$.  Hence $\Gamma_i\subset\Gamma_j$ and $\partial A_{2i+1}\subset\Gamma_j$.  Note that $\Gamma_j\subset S$ is an annulus carried by $N(B)$ and $\partial A_{2i+1}\subset\Gamma_j\cap A$, so a sub-annulus of $\Gamma_j$ satisfies the hypotheses of Proposition~\ref{Ptorus}.  Hence $B$ must carry a torus, contradicting our hypotheses.  Thus, the solid tori $T_i$'s are pairwise disjoint.  Note that $\partial T_i\subset N(B)$ but the solid torus $T_i$ is not contained in $N(B)$, since $A_k\subset A\subset N(B)$ is a vertical annulus.  So each solid torus $T_i$ must contain a component of $\partial_hN(B)$, and hence the number of such solid tori is bounded by the number of components of $\partial_hN(B)$.  Therefore, there is a number $K$ depending only on $B$ such that $n=|A\cap S|<K$.

If $c_1$ bounds a compressing disk in $H_1$, since $c_1\cup c_2$ bounds the annulus $A_1$ in $H_1$, $c_2$ bounds a compressing disk in $H_1$.  As $A_2$ is an annulus in $H_2$, by Lemma~\ref{Lnest}, we have that $A_2$ is $\partial$--parallel in $H_2$.  Using the same argument, we can inductively conclude that $A_{2k}$ is $\partial$--parallel in $H_2$ for each $k$, and obtain such a bound $K$ on $n=|A\cap S|$.
\end{proof}

The following Proposition for branched surfaces is well-known,  see also \cite{FO,AL}.

\begin{proposition}\label{Pmon}
Let $B$ be a branched surface in $M$.  Suppose $M-B$ is irreducible and $\partial_hN(B)$ is incompressible in $M-int(N(B))$.  Let $C$ be a component of $M-int(N(B))$ and suppose $C$ contains a monogon.  Then $C$ must be a solid torus
in the form of $D\times S^1$, where $D$ is a monogon.
\end{proposition}
\begin{proof}
Let $D$ be a monogon in $C$, {\it i.e.}, the disk $D$ is properly embedded in $C$,  $\partial D$ consists of two arcs, $\alpha\subset\partial_vN(B)$ and $\beta\subset\partial_hN(B)$, and $\alpha$ is a vertical arc in $\partial_vN(B)$.  Let $v$ be the component of $\partial_vN(B)$ containing $\alpha$.  Then as shown in Figure~\ref{monogon} (a), the union of two parallel copies of $D$ and a rectangle in $v$ is a disk $E$ properly embedded in $C$, with $\partial E\subset \partial_hN(B)$.  Since $\partial_hN(B)$ is incompressible in $M-int(N(B))$, $\partial E$ must bound a disk in $\partial_hN(B)\cap\partial C$.  Since $C$ is irreducible, $C$ must be a solid torus in the form of $D\times S^1$, where $D$ is the monogon above.
\end{proof}

Before we proceed, we quote two results of Scharlemann that we will use later.

\begin{lemma}[Lemma 2.2 of \cite{S}]\label{Ls22}
Suppose $H_1\cup_SH_2$ is a strongly irreducible Heegaard splitting of a 3--manifold $M$ and $F$ is a disk in $M$ transverse to $S$ with $\partial F\subset S$.  Then $\partial F$ bounds a disk in some $H_i$.
\end{lemma}

\begin{theorem}[Theorem 2.1 of \cite{S}]\label{Ts21}
Suppose $H_1\cup_SH_2$ is a strongly irreducible Heegaard splitting of a 3--manifold $M$ and $B$ is 3--ball in $M$.  Let $T_i$ be the planar surface $\partial B\cap H_i$ properly embedded in $H_i$, and suppose $T_i$ is incompressible in $H_i$.  Then $S\cap B$ is connected and $\partial$--parallel in $B$.
\end{theorem}

Corollary~\ref{Csch} follows trivially from Scharlemann's theorem.

\begin{corollary}\label{Csch}
Suppose $H_1\cup_SH_2$ is a strongly irreducible Heegaard splitting of a 3--manifold.  Let $P$ be a planar surface properly embedded in $H_1$.  Suppose $P$ is incompressible in $H_1$, and each boundary component of $P$ bounds a disk in $H_2$.  Then $P$ is $\partial$--parallel in $H_1$.
\end{corollary}

\section{Measured laminations}\label{Slam}

The purpose of this section is to prove Lemma~\ref{Lfine}, which is an easy consequence of some properties of laminations and results from \cite{L4}.

The following theorem is one of the fundamental results in the theory of measured laminations and foliations. It also plays an important role in \cite{L4}.  An exceptional minimal lamination is a lamination in which every leaf is dense, and the intersection of any transversal with such a lamination is a Cantor set, see section 3 of \cite{L4}.

\begin{theorem}[Theorem 3.2 in Chapter I of \cite{MS}, pp 410]\label{TMS}
Let $\mu$ be a co-dimension one measured lamination in a closed connected 3--manifold $M$, and suppose $\mu\ne M$.  Then $\mu$ is the disjoint union of a finite number of sub-laminations.  Each of these sub-laminations is of one of the following types:
\begin{enumerate}
\item A family of parallel compact leaves,
\item A twisted family of compact leaves,
\item An exceptional minimal measured lamination.
\end{enumerate} 
\end{theorem}

\begin{definition}[Definition 4.2 of \cite{L4}]\label{Dvan}
Let $\mu$ be a lamination in $M$ and $l_0$ a leaf of $\mu$.  We call a simple closed curve $f_0:S^1\to l_0$ an \emph{embedded vanishing cycle} in $\mu$ if $f_0$ extends to an embedding $F:[0, 1]\times S^1\to M$ satisfying the following properties.
\begin{enumerate}
\item $F^{-1}(\mu)=C\times S^1$, where $C$ is a closed set of $[0, 1]$, and for any $t\in C$, the curve $f_t(S^1)$, defined by $f_t(x)=F(t,x)$, is contained in a leaf $l_t$, 
\item for any $x\in S^1$, the curve $t\to F(t,x)$ is transverse to $\mu$,
\item $f_0$ is an essential curve in $l_0$, but there is a sequence of points $\{t_n\}$ in $C$ such that $\lim_{n\to\infty}t_n=0$ and $f_{t_n}(S^1)$ bounds a disk in $l_{t_n}$ for all $t_n$. 
\end{enumerate}
\end{definition}

The following lemma from \cite{L4} will be useful in our proof of Lemma~\ref{Lfine}.

\begin{lemma}[Lemma 4.3 of \cite{L4}]\label{Lvan}
Let $M$ be a closed orientable and irreducible 3--manifold, and $\mu\subset M$ an exceptional minimal measured lamination.  Suppose $\mu$ is fully carried by a branched surface $B$ and $B$ does not carry any $2$--sphere.  Then $\mu$ has no embedded vanishing cycle.
\end{lemma}

The proof of the follow lemma is similar in spirit to part of the proof of Lemma 4.5 in \cite{L4}.

\begin{lemma}\label{Lleaf}
Let $B$ be a branched surface in a closed, orientable and irreducible 3--manifold $M$, and $M\ne T^3$.  Suppose $B$ does not carry any 2--sphere or torus, and suppose $B$ fully carries a measured lamination $\mu$.  Then $\mu$ does not contain any plane leaf, infinite annular leaf or infinite M\"{o}bius band leaf.
\end{lemma}
\begin{proof}
By Theorem~\ref{TMS}, we may assume $\mu$ is an exceptional minimal measured lamination, in particular, every leaf is dense in $\mu$.

Suppose every leaf of $\mu$ is a plane.  After trivially eliminating all the disks of contacts in $N(B)$ that are disjoint from $\mu$, we have that $\partial_hN(B)$ consists of disks.  So there is no monogon and $\mu$ is an essential lamination.  By a Theorem in \cite{G7} (also see Proposition 4.2 of \cite{L1}), $M\cong T^3$.

So at least one leaf of $\mu$ is not a plane.  Let $\gamma$ be an essential simple closed curve in a non-plane leaf.  Since $\mu$ is a measured lamination, there is no holonomy.  So there is an embedded vertical annulus $S^1\times I\subset N(B)$ such that $\gamma\subset S^1\times I$ and $\mu\cap(S^1\times I)$ is a union of parallel circles.  Suppose $L$ is a plane leaf of $\mu$. Since every leaf is dense, $L\cap (S^1\times I)$ contains infinitely many circles whose limit is $\gamma$. As $L$ is a plane, these circles bound disks in $L$.  By Definition~\ref{Dvan}, $\gamma$ is an embedded vanishing cycle, and this contradicts Lemma~\ref{Lvan}.  So $\mu$ does not contain any plane leaf.

Suppose $\mu\subset N(B)$ and $A$ is an infinite annular leaf (or an infinite M\"{o}bius band leaf) of $\mu$.   Let $\gamma$ be an essential simple closed curve in $A$.  There is an embedded vertical annulus $S^1\times I\subset N(B)$ such that $\gamma\subset S^1\times I$, and $\mu\cap (S^1\times I)$ is a union of parallel circles.  Since every leaf is dense in $\mu$, $A\cap (S^1\times I)$ contains infinitely many circles whose limit is $\gamma$.  By Lemma~\ref{Lvan}, we may assume that only finitely many circles of $A\cap (S^1\times I)$ are trivial in $A$.  So there exist 3 essential simple closed curves in $A\cap (S^1\times I)$, $\gamma_i$ ($i=1,2,3$), such that $\gamma_1\cup\gamma_3$ bounds a compact sub-annulus $A_\gamma$ in $A$ with $int(A_\gamma)\cap (S^1\times I)=\gamma_2$.  By Proposition~\ref{Ptorus}, $B$ carries a torus, contradicting our hypotheses.
\end{proof}

\begin{lemma}\label{Linject}
Let $B$ be a branched surface in $M$.  Suppose $N(B)$ does not contain any disk of contact and $\partial_hN(B)$ has no disk component.  Let $\lambda\subset N(B)$ be a lamination fully carried by $N(B)$.  Then every leaf of $\lambda$ is $\pi_1$--injective in the 3--manifold $N(B)$.
\end{lemma}
\begin{proof}
We may use the arguments in \cite{GO} to prove this lemma directly, but it is more convenient to simply use a theorem of \cite{GO}.  Since $\partial_hN(B)$ has no disk component, no component of $\partial N(B)$ is a 2--sphere.  For each component $S$ of $\partial N(B)$, we may glue to $N(B)$ (along $S$) a compact orientable and irreducible 3--manifold $M_S$, whose boundary $\partial M_S\cong S$ is incompressible in $M_S$.  So we can obtain a closed 3--manifold $M'$ this way with $N(B)\subset M'$.  Since $S$ is $\pi_1$--injective in $M_S$, the inclusion $i:N(B)\hookrightarrow M'$ induces an injection on $\pi_1$.

If $\partial_hN(B)$ is compressible in $M'-int(N(B))$, then we have a compressing disk $D$ with $\partial D\subset\partial_hN(B)\cap S$, where $S$ is a boundary component of $N(B)$.  As $S$ is incompressible in $M_S$, $\partial D$ must bound a disk $E$ in $S$, which implies that $E$ contains a disk component of $\partial_hN(B)$, contradicting our hypotheses.  So $\partial_hN(B)$ must be incompressible in $M'-int(N(B))$.  There is clearly no monogon by the construction and no disk of contact by our hypotheses.  Moreover, since $\partial_hN(B)$ has no disk component and there is no monogon, it is easy to see that there is no Reeb component for $N(B)$.  Therefore, by \cite{GO}, $\lambda$ is an essential lamination in the closed manifold $M'$, and every leaf of $\lambda$ is $\pi_1$--injective in $M'$ hence $\pi_1$--injective in $N(B)$.
\end{proof}

The following lemma from \cite{L4} is also useful in the proof of Lemma~\ref{Lfine}.  

\begin{lemma}[Lemma 4.1 of \cite{L4}]\label{Lnodoc}
Let $B$ be a branched surface fully carrying a lamination $\mu$.  Suppose $\partial_hN(B)$ has no disk component and $N(B)$ does not contain any disk of contact that is disjoint from $\mu$.  Then $N(B)$ does not contain any disk of contact.
\end{lemma}

Now, Lemma~\ref{Lfine} follows easily from the previous lemmas.

\begin{lemma}\label{Lfine}
Let $B$ be a branched surface in a closed, orientable and irreducible 3--manifold $M$.  Suppose $B$ does not carry any 2--sphere or torus, and $B$ fully carries a measured lamination $\mu$.  Then $B$ can be split into a branched surface $B_1$ such that $B_1$ still fully carries $\mu$, no component of $\partial_hN(B_1)$ is a disk, and every leaf of $\mu$ is $\pi_1$--injective in $N(B_1)$.
\end{lemma}
\begin{proof}
By Theorem~\ref{TMS}, we may assume that $\mu$ is an exceptional minimal measured lamination.  Since $B$ does not carry any 2--sphere or torus, by Lemma~\ref{Lleaf}, no leaf of $\mu$ is a plane.  After some isotopy, we may assume $\partial_hN(B)\subset\mu$.  Hence we can split $N(B)$ so that each component of $\partial_hN(B)$ contains an essential curve of the corresponding leaf.  So no component of $\partial_hN(B)$ is a disk after the splitting.

By splitting $N(B)$, we may trivially eliminate all the disks of contact that are disjoint from $\mu$.  So, by Lemma~\ref{Lnodoc}, $N(B)$ does not contain any disk of contact.  Now the lemma follows from Lemma~\ref{Linject}.
\end{proof}

The following Proposition is well-known.  It also plays a fundamental role in \cite{L5}.

\begin{proposition}\label{PHaken}
Let $M$ be a closed irreducible and orientable 3--manifold and $B$ a branched surface in $M$ carrying  a measured lamination $\mu$.  If $\mu$ is an essential lamination, then $B$ carries an incompressible surface and hence $M$ is Haken.
\end{proposition}
\begin{proof}
By \cite{GO}, if $\mu$ is an essential lamination, then one can split $B$ into an incompressible branched surface $B'$ that fully carries $\mu$.  Since $\mu$ is a measured lamination, the system of branch equations for $B'$ must have a positive solution.  Since the coefficients of each branch equation are integers, the system of branch equations must have a positive integer solution.  Thus $B'$ fully carries a closed orientable surface.  By \cite{FO}, every closed surface fully carried by an incompressible branched surface is incompressible.
\end{proof}

\section{Limits of compact surfaces}\label{Slimit}

Let $B$ be a branched surface in a closed 3--manifold $M$, and $F\subset N(B)$ a closed surface carried by $B$.  Then $F$ corresponds to a non-negative integer solution to the branch equations of $B$, see section 3 of \cite{L4} for a brief explanation and see \cite{FO,O} for more details.  We use $\mathcal{S}(B)\subset\mathbb{R}^N$ to denote the set of non-negative solutions to the branch equations of $B$, where $N$ is the number of branch sectors of $B$.  There is a one-to-one correspondence between a closed surface carried by $B$ and an integer point in $\mathcal{S}(B)$.  A surface is fully carried by $B$ if and only if every coordinate of the corresponding point in $\mathcal{S}(B)$ is positive.  

Every point in $\mathcal{S}(B)$, integer point or non-integer point, corresponds to a measured lamination carried by $B$.  Such a measured lamination $\mu$ can be viewed as the inverse limit of a sequence of splittings $\{B_n\}_{n=0}^\infty$, where $B_0=B$ and $B_{i+1}$ is obtained by splitting $B_i$. Note that if $B_{i+1}$ is obtained by splitting $B_i$, one may naturally consider $N(B_{i+1})\subset N(B_i)$.  We refer  to section 3 of \cite{L4} for a brief description, see \cite{O} and section 3 of \cite{Hat} for more details (also see Definition 4.1 and Lemma 4.2 of \cite{GO}).  There is a one-to-one correspondence between each point in $\mathcal{S}(B)$ and a measured lamination constructed in this fashion.  This one-to-one correspondence is slightly different from the one above for integer points of $\mathcal{S}(B)$.  For an integer point, the sequence of splittings on $B$ above stop in a finite number of steps (i.e., $B_{i+1}=B_i$ is a closed surface if $i$ is large), and the measured lamination constructed this way is the horizontal foliation of an $I$--bundle over a closed surface. 

We define the \emph{projective lamination space} of $B$, denoted by $\mathcal{PL}(B)$, to be the set of points in $\mathcal{S}(B)$ satisfying $\sum_{i=1}^Nx_i=1$.  Let $p: \mathcal{S}(B)-\{0\}\to\mathcal{PL}(B)$ be the natural projection sending $(x_1,\dots,x_N)$ to $\frac{1}{s}(x_1,\dots,x_N)$, where $s=\sum_{i=1}^Nx_i$.  To simplify notation, we do not distinguish a point $x\in\mathcal{S}(B)$ and its image $p(x)\in\mathcal{PL}(B)$ unless necessary. $\mathcal{PL}(B)$ is a compact set. For any infinite sequence of distinct closed surfaces carried by $B$, the images of the corresponding points in $\mathcal{PL}(B)$ (under the map $p$) has an accumulation point, which corresponds to a  measured lamination $\mu$.  To simplify notation, we simply say that the measured lamination $\mu$ is an accumulation point of this sequence of surfaces in $\mathcal{PL}(B)$.  Throughout this paper, when we consider a compact surface carried by $B$, we identify the surface with an integer point in $\mathcal{S}(B)$, but when we consider $\mu$ as a limit point of a sequence of compact surfaces in $\mathcal{PL}(B)$, we identify the point $\mu\in\mathcal{PL}(B)$ to a measured lamination as the inverse limit of the sequence of splittings on $B$ above.

\begin{proposition}\label{Plinear}
Let $B$ be a branched surface with $n$ branch sectors and $\{S_k=(x_1^{(k)}, \dots, x_n^{(k)})\}$ an infinite sequence of integer points in $\mathcal{S}(B)$ whose images in $\mathcal{PL}(B)$ are distinct points.  Suppose $\mu=(z_1,\dots,z_n)\in\mathcal{PL}(B)$ is the limit point of $\{S_k\}$ in the projective lamination space.  Let $f(x_1,\dots,x_n)$ be a homogeneous linear function with $n$ variables.  Then we have the following.
\begin{enumerate}
\item If $z_i=0$ and $z_j\ne 0$, then $\lim_{k\to\infty}x_i^{(k)}/x_j^{(k)}=0$.
\item If $z_i>z_j$, then $x_i^{(k)}>x_j^{(k)}$ if $k$ is sufficiently large.
\item If the sequence $\{f(S_k)\}$ is bounded, then $f(\mu)=0$.
\end{enumerate}
\end{proposition}
\begin{proof}
Let $s_k=\sum_{i=1}^nx_i^{(k)}$.  Then the corresponding point of $S_k$ in $\mathcal{PL}(B)$ is $[S_k]=(x_1^{(k)}/s_k,\dots,x_n^{(k)}/s_k)$. By our hypotheses, $\lim_{k\to\infty}x_i^{(k)}/s_k=z_i$ for each $i$.  Thus, if $z_i=0$ and $z_j\ne 0$, we have $\lim_{k\to\infty}x_i^{(k)}/x_j^{(k)}=z_i/z_j=0$.

Since  $x_i^{(k)}/s_k>x_j^{(k)}/s_k$ is equivalent to $x_i^{(k)}>x_j^{(k)}$, part 2 is obvious.

Since $f(x_1,\dots,x_n)$ is a homogeneous linear function, $f([S_k])=f(S_k)/s_k$ and $\lim_{k\to\infty}f([S_k])=f(\mu)$.  Since the sequence $\{S_k=(x_1^{(k)},\dots,x_n^{(k)})\}$ consists of distinct non-negative integer solutions, the integers $\{s_k\}$ are unbounded. So, after passing to a sub-sequence if necessary, we have $\lim_{k\to\infty}s_k=\infty$.  Therefore, if the sequence $\{f(S_k)\}$ is bounded from above, then $\lim_{k\to\infty}f(S_k)/s_k=f(\mu)=0$.
\end{proof}

\begin{corollary}\label{Clinear}
Let $\{S_k\}\subset N(B)$ be a sequence of distinct compact connected surfaces carried by a branched surface $B$.  Suppose $\mu\subset N(B)$ is the measured lamination corresponding to the limit of $\{S_k\}$ in $\mathcal{PL}(B)$, and let $K$ be an $I$--fiber of $N(B)$ such that $K\cap\mu\ne\emptyset$.  Then, if $k$ is large, $|K\cap S_k|$, the number of intersection points of $K$ and $S_k$, is large.
\end{corollary}
\begin{proof}
The number of intersection points of an $I$--fiber and $S_k$ is equal to the integer value of a coordinate of the corresponding point in $\mathcal{S}(B)$.  So the corollary follows immediately from part 3 of Proposition~\ref{Plinear} after setting the linear function to $f(x_1,\dots,x_n)=x_i$, where $x_i$ corresponds to the branch sector of $B$ that contains the point $\pi(K)$ ($x_i=|K\cap S_k|$).
\end{proof}

We call a lamination $\mu$ a normal lamination with respect to a triangulation if every leaf of $\mu$ is a (possibly non-compact) normal surface.

\begin{corollary}\label{Cnormal}
Let $M$ be a closed 3--manifold with a fixed triangulation, and let $B$ be a branched surface obtained by gluing together a collection of normal disks and at most one almost normal piece, similar to \cite{FO}.  Suppose $\{S_n\}$ is an infinite sequence of distinct connected almost normal surfaces fully carried by $B$. Then each accumulation point of $\{S_n\}$ in $\mathcal{PL}(B)$ must correspond to a normal measured lamination.
\end{corollary}
\begin{proof}
If $B$ does not contain an almost normal piece, then every surface carried by $B$ is normal and there is nothing to prove.  Suppose $s$ is a branch sector of $B$ containing the almost normal piece.  Since $B$ fully carries an almost normal surface, $B-int(s)$ must be a sub-branched surface of $B$ and every lamination carried by $B-int(s)$ is normal ($B-int(s)$ is called the normal part of $B$ in section 2 of \cite{L4}).  Suppose $S_n=(x_1,\dots,x_N)\in\mathcal{S}(B)$ and suppose $x_1$ is the coordinate corresponding to the branch sector $s$.  Since an almost normal surface has at most one almost normal piece, $x_1=1$ for each $S_n$.  Suppose $\mu=(z_1,\dots,z_N)\in\mathcal{PL}(B)$.  By Proposition~\ref{Plinear} and Corollary~\ref{Clinear}, $z_1$ must be zero.  Hence $\mu$ is carried by $B-int(s)$ and is a normal lamination.
\end{proof}

Now, we will use two examples to illustrate the limit of closed surfaces.  Although the two examples are train tracks, similar results hold for branched surfaces.

\begin{figure}
\begin{center}
\includegraphics{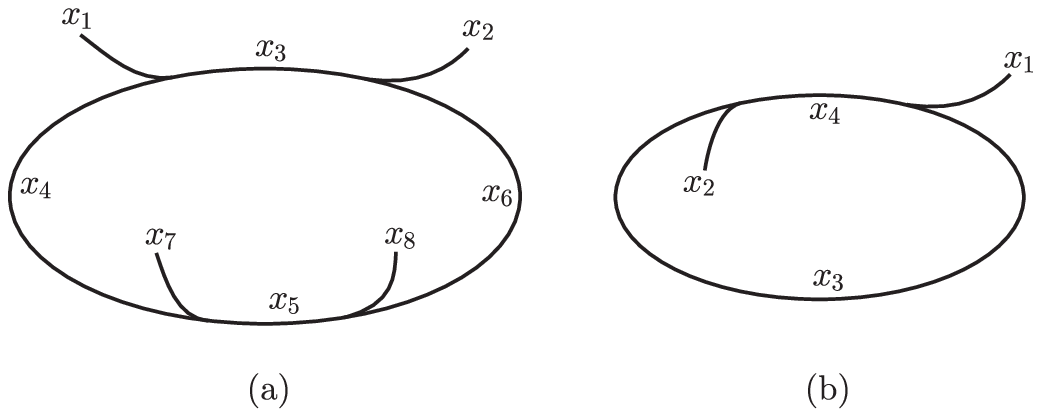}
\caption{}\label{traintrack}
\end{center}
\end{figure}

\begin{example}\label{Extrain}
Let $\tau$ be a train track in the plane as shown in Figure~\ref{traintrack}(a).  There are 8 branch sectors in $\tau$, and the branch equations are $x_1+x_4=x_3=x_2+x_6$ and $x_7+x_4=x_5=x_8+x_6$.  Suppose $\{\gamma_n\}$ is an infinite sequence of compact arcs carried by $\tau$ whose limit in $\mathcal{PL}(\tau)$ is the point $\mu=(0,0,1/4,1/4,1/4,1/4,0,0)$.  Geometrically $\mu$ is a measured lamination consisting of parallel circles carried by $\tau$. We identify $\gamma_n$ with an integer point in $\mathcal{S}(\tau)$ and suppose the $\gamma_n$'s are different points in $\mathcal{S}(\tau)$. Note that $\gamma_n$ contains a circle if and only if $x_1=x_2$ and $x_7=x_8$.  By Proposition~\ref{Plinear} and Corollary~\ref{Clinear}, as $n$ tends to infinity, the values $x_6$ and $x_6/x_2$ of $\gamma_n$ tend to infinity.  This implies that, if $n$ is large, $\gamma_n$ contains either many parallel circles or a spiral wrapping around the circle many times.
\end{example}

In the proof of the main theorem, we will consider the limit lamination $\mu$ of an infinite sequence of almost normal Heegaard surfaces carried by a branched surface $B$.  The measured lamination $\mu$ is fully carried by a sub-branched surface $B^-$ of $B$.   In many situations, we would like to split $B^-$ into a nicer branched surface $B^-_1$. In fact, by considering $\mu\subset N(B^-)\subset N(B)$, we can split $N(B^-)$ and $N(B)$ simultaneously and obtain $\mu\subset N(B_1^-)\subset N(B_1)$, such that $B_1^-$ is the sub-branched surface of $B_1$ that fully carries $\mu$, $B_1$ is obtained by splitting $B$, and $B_1$ still carries an infinite sub-sequence of $\{S_n\}$.  Next, we will use Example~\ref{Exsplit} to illustrate how the local splittings work.  We also formulate this fact in Proposition~\ref{Psplit}.  Proposition~\ref{Psplit} is similar in spirit to Lemma 6.1 of \cite{L4}.

\begin{figure}
\begin{center}
\includegraphics{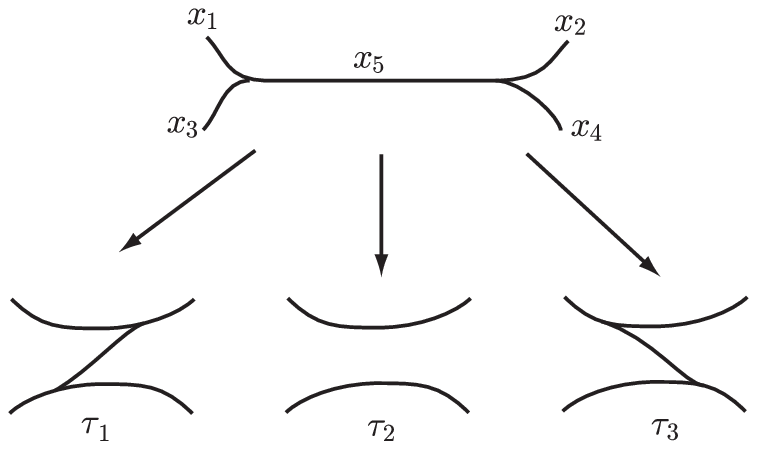}
\caption{}\label{splitting}
\end{center}
\end{figure}

\begin{example}\label{Exsplit} 
Let $\tau$ be the train track on the top of Figure~\ref{splitting}.  As shown in Figure~\ref{splitting}, $\tau$ can be split into 3 different train tracks $\tau_1$, $\tau_2$ and $\tau_3$.  Suppose $\mu$ is a lamination fully carried by $\tau$. Let $x_1,\dots,x_5$ be the weights of $\mu$ at the branch sectors of $\tau$.  These $x_i$'s satisfy the branch equations $x_1+x_3=x_5=x_2+x_4$.  It is easy to see that $x_1<x_2$ (resp. $x_1>x_2$) if and only if $\mu$ is fully carried by $\tau_1$ (resp. $\tau_3$), and $x_1=x_2$ if and only if $\mu$ is fully carried by $\tau_2$.  Suppose $\{S_n\}$ is an infinite sequence of compact arcs carried by $\tau$ and suppose each $S_n$ corresponds to a different integer point in $\mathcal{S}(\tau)$.  Suppose the limit of $\{S_n\}$ in $\mathcal{PL}(\tau)$ is $\mu$.  

By part 2 of Proposition~\ref{Plinear}, if $x_1<x_2$ (resp. $x_1>x_2$) for $\mu$, we can split $\tau$ into $\tau_1$ (resp. $\tau_3$), and $\tau_1$ (resp. $\tau_3$) fully carries $\mu$ and an infinite sub-sequence of $\{S_n\}$.  Now, we consider the case $x_1=x_2$ for $\mu$.  Although we can split $\tau$ (along $\mu$) into $\tau_2$ which fully carries $\mu$, $\tau_2$ may not carry infinitely many $S_n$'s. Nonetheless, if $\tau_2$ only carries finitely many $S_n$'s, then at least one of $\tau_1$ and $\tau_2$, say $\tau_1$, must carry an infinite sub-sequence of $\{S_n\}$.  Moreover, $\tau_1$ can be considered as the train track obtained by adding a branch sector to $\tau_2$, and $\tau_1$ can be obtained by splitting $\tau$.  

Now we consider the splittings on branched surfaces.  Note that any splitting on a branched surface can be viewed as a sequence of successive local splittings, and the operations of such local splittings on a branched surface are basically the same as the splittings on the train track in Example~\ref{Exsplit}.  So we have the following proposition.
\end{example}

\begin{proposition}\label{Psplit}
Let $B$ be a branched surface and $\{S_n\}\subset\mathcal{S}(B)$ a sequence of distinct compact surfaces carried by $B$.  Suppose $\mu\in\mathcal{PL}(B)$ is the limit point of $\{S_n\}$ in $\mathcal{PL}(B)$.  Let $B^-$ be the sub-branched surface of $B$ that fully carries $\mu$.  Let $B^-_1$ be any branched surface obtained by splitting $B^-$ along $\mu$, and suppose $B^-_1$ still fully carries $\mu$.  Then one can add some branch sectors to $B^-_1$ to form a branched surface $B_1$ (i.e. $B_1^-$ is a sub-branched surface of $B_1$), such that $B_1$ can be obtained by splitting $B$, and $B_1$ carries an infinite sub-sequence of $\{S_n\}$.
\end{proposition}
\begin{proof}
This proposition is similar in spirit to Lemma 6.1 of \cite{L4}.  The splitting from $B^-$ to $B_1^-$ can be divided into a sequence of successive small local splittings, and each local splitting is similar to the splittings in Example~\ref{Exsplit} and Figure~\ref{splitting}.  During each local splitting, we can split $B^-$ and $B$ simultaneously.  If $B$ fails to carry infinitely many $S_n$'s after a local splitting, similar to the operation of obtaining $\tau_1$ by adding a branch sector to $\tau_2$ in Example~\ref{Exsplit}, we can always add some branched sectors to get a branched surface satisfying the requirements of the Proposition.
\end{proof}
\begin{remark}
In Proposition~\ref{Psplit}, $B_1^-$ is the sub-branched surface of $B_1$ that fully carries $\mu$.  Since any lamination carried by $B_1$ is carried by $B$, it is easy to see that $\mu\subset N(B_1)$ is still the limit point in $\mathcal{PL}(B_1)$ of the sub-sequence of $\{S_n\}$ carried by $B_1$. 
\end{remark}

\begin{remark}\label{Rsss}
Let $\{S_n\}$ be an infinite sequence of distinct closed surfaces carried by $N(B)$ whose limit in $\mathcal{PL}(B)$ is a measured lamination $\mu$.  Let $\gamma$ be a simple closed essential curve in a leaf of $\mu$.  If every $I$--fiber of $N(B)$ intersects $\gamma$ in at most one point, then (after a slight enlargement) $\pi^{-1}(\pi(\gamma))$ can be considered as a fibered neighborhood of a train track consisting of a circle $\pi(\gamma)$ and some ``tails" along the circle similar to Figure~\ref{traintrack}, where $\pi:N(B)\to B$ is the map collapsing each $I$--fiber to a point.  Since the limit of $\{S_n\}$ is $\mu$, $\pi^{-1}(\pi(\gamma))\cap S_i$ ($i=1,2,\cdots$) is a sequence of curves whose limit is a measured lamination containing the circle $\gamma$.  As in Example~\ref{Extrain}, if $n$ is large, $\pi^{-1}(\pi(\gamma))\cap S_n$ must contain either many circles parallel to $\gamma$ or a spiral winding around $\gamma$ many times.  However, if there is an $I$--fiber of $N(B)$ intersecting $\gamma$ in more than one point, then $\pi(\gamma)$ is an immersed curve in $B$.  Nevertheless, since $\gamma$ is an embedded essential curve in a leaf of $\mu$, by Theorem~\ref{TMS}, after a finite sequence of splittings on $B$, we can get a branched surface $B_1$ such that $B_1$ still carries $\mu$ and $\pi|_\gamma$ is injective, where $\pi:N(B_1)\to B_1$ is the collapsing map, (i.e., every $I$--fiber of $N(B_1)$ intersects $\gamma$ in at most one point).  Moreover, by Proposition~\ref{Psplit}, we may assume $B_1$ still carries an infinite sub-sequence of $\{S_n\}$.  Now the situation is the same as above after replacing $B$ by $B_1$.
\end{remark}

The next lemma says that, if the branched surface is nice, then the limit of trivial circles in a sequence of closed surfaces cannot be an essential circle in the limit lamination.

\begin{lemma}\label{Lessential}
Let $M$ be a closed 3--manifold with a fixed triangulation, and let $B$ be a branched surface obtained by gluing together a collection of normal disks and at most one almost normal piece, as in Theorem~\ref{Heeg1}. Suppose $N(B)$ does not carry any normal 2--sphere or normal torus.  Let $\{S_n\}$ be a sequence of distinct closed almost normal surfaces fully carried by $N(B)$ whose limit in $\mathcal{PL}(B)$ is a measured lamination $\mu$.   Let $\gamma$ be an essential simple closed curve in a leaf of $\mu$.  Then $B$ can be split into a branched surface $B_1$ that carries both $\mu$ and an infinite sub-sequence $\{S_{n_k}\}$ of $\{S_n\}$, such that, for any embedded vertical annulus $A\supset\gamma$ in $N(B_1)$, $A\cap S_{n_k}$ does not contain any circle that is trivial in the surface $S_{n_k}$, for each $S_{n_k}$.
\end{lemma}
\begin{proof}
Let $A_\gamma$ be an embedded vertical annulus in $N(B)$ containing $\gamma$.  Suppose $A_\gamma\cap S_n$ contains a trivial circle in $S_n$ for each $n$.  Such a trivial circle bounds a disk $D_n$ in $S_n$.  So $D_n$ is transverse to the $I$--fibers of $N(B)$ and with $\partial D_n\subset A_\gamma$.  Let $s$ be the branch sector containing the almost normal piece, and let $B'=B-int(s)$ be the sub-branched surface of $B$ ($B'$ is called the normal part of $B$, see section 2 of \cite{L4}).  By Corollary~\ref{Cnormal}, $\mu$ is carried by $B'$.  So we can assume that if $D_n$ contains an almost normal piece, the almost normal piece lies in $int(D_n)$.  Since $S_n$ is an almost normal surface, $D_n$ contains at most one almost normal piece.  

We call an isotopy of $N(B)$ a $B$--isotopy, if the isotopy is invariant on each $I$--fiber of $N(B)$.

\medskip
\noindent\textbf{Claim}. Up to $B$--isotopy, there are only finitely many such disks $D_n$.
\medskip

To prove the claim, we first consider such disks that do not contain almost normal pieces.  If $D_n$ does not contain an almost normal piece, then we may assume that $D_n$ lies in $N(B')$ transverse to the $I$--fibers of $N(B')$, and consider $A_\gamma$ as an embedded vertical annulus in $N(B')$.  Let $\mathcal{S}_\gamma$ be the set of embedded compact surfaces $F$ in $N(B')$ with the properties that $F$ is transverse to the $I$--fibers of $N(B')$ and $\partial F$ is a single circle in $A_\gamma$.  Similar to $\mathcal{S}(B')$, we can describe $\mathcal{S}_\gamma$ as the set of non-negative integer solutions of a system of non-homogeneous linear equations as follows, see \cite{AL} for such a description for disks of contact.  Let $L'$ be the branch locus of $B'$ and suppose $\pi(A_\gamma)$ is an immersed curve in $B'$.  Suppose $b_1,\dots, b_N$ are the components of $B'-L'-\pi(A_\gamma)$.  For each $b_i$ and any $F\in\mathcal{S}_\gamma$, let $x_i=|F\cap\pi^{-1}(b_i)|$.  One can describe $F$ using a non-negative integer point $(x_1,\dots,x_N)\in\mathbb{R}^N$, and  $(x_1,\dots,x_N)$ is a solution of the system of (non-homogeneous) linear equations in the forms of $x_k=x_i+x_j$ and $x_i=x_j+1$.  Equations like $x_i=x_j+1$ occur when two pieces are glued along $\pi(A_\gamma)$, since $\pi(\partial F)=\pi(A_\gamma)$. Up to $B'$--isotopy, there are only finitely many surfaces corresponding to the same integer point in $\mathcal{S}_\gamma$.  Moreover, the corresponding homogeneous system is exactly the system of branch equations of $B'$.  Suppose there is an infinite sequence of distinct disks $\{D_n\}$ in $\mathcal{S}_\gamma$. Then one can find $D_i=(x_1,\dots,x_N)$ and $D_j=(y_1,\dots,y_N)$ such that $x_k\le y_k$ for each $k$.  Thus $D_j-D_i=(y_1-x_1,\dots,y_N-x_N)$ is a non-negative integer solution to the corresponding homogeneous system, i.e., the system of branch equations.  So $D_j-D_i$ corresponds to a closed surface carried by $B'$.  Since the Euler characteristic is additive, $\chi(D_j-D_i)=\chi(D_j)-\chi(D_i)=0$. This means $B'$ carries a closed surface (which may not be connected) with total Euler characteristic $0$, which implies that $B'$ must carry a connected surface with non-negative Euler characteristic.  If $B'$ carries a Klein bottle (or projective plane), $B'$ must carry a torus (or $2$--sphere) because $M$ is orientable.  Since $B'=B-int(s)$, every surface carried by $B'$ is normal.  This contradicts the hypothesis that $B$ does not carry any normal 2--sphere or normal torus. 

Suppose there is an infinite sequence of disks $\{D_n\}$ from the $S_n$'s, such that each $D_n$ contains an almost normal piece.  As above, we can also identify each $D_n$ as an integer solution of a system of non-homogeneous linear equations. Up to $B$--isotopy, there are only finitely many such disks corresponding to the same integer point.  If the disks $\{D_n\}$ correspond to different integer points, then one can find $D_i=(x_1,\dots,x_K)$ and $D_j=(y_1,\dots,y_K)$ such that $x_k\le y_k$ for each $k$.  Suppose the first coordinate corresponds to the branch sector $s$ that contains the almost normal piece.  Since each $S_n$ is an almost normal surface, each $D_n$ contains only one almost normal piece.  Hence, $x_1=y_1=1$ and the first coordinate of $D_j-D_i$ is $y_1-x_1=0$.  This means that $D_j-D_i$ does not contain an almost normal piece and is carried by $B'$.  Now the argument is the same as above.  This finishes the proof of the claim.

Let $B^-$ be the sub-branched surface of $B$ fully carrying $\mu$.  As described earlier in this section and in section 3 of \cite{L4}, we may consider $\mu$ as the inverse limit of an infinite sequence of splittings on $B^-$.  Suppose $\{B_n^-\}_{n=0}^\infty$ ($B_0^-=B^-$) is such a sequence of branched surfaces, with each $B_i^-$ obtained by splitting $B_{i-1}^-$ and $\mu$ being the inverse limit of the sequence $\{N(B_n^-)\}$.  Note that if $\mu$ consists of compact leaves, then such splittings are a finite process.  By Theorem~\ref{TMS}, we only consider the case that $\mu$ is an exceptional minimal measured lamination, and the proof for the case that $\mu$ consists of compact leaves is similar.  By Proposition~\ref{Psplit}, we may assume there is a sequence of branched surfaces $\{B_n\}$ ($B_0=B$) such that, for each $n$, $B_{n+1}$ is obtained by splitting $B_n$, $B_n$ carries $\mu$ and an infinite sub-sequence of $\{S_n\}$, and $B_n^-$ is a sub-branched surface of $B_n$.

Let $A_k\subset N(B_k^-)$ be a vertical annulus containing $\gamma$.  By Lemma~\ref{Lfine}, after some splittings, we may assume that if $k$ is sufficiently large, every leaf of $\mu$ is $\pi_1$--injective in $N(B_k^-)$.  Since $\gamma$ is an essential curve in a leaf, if $k$ is sufficiently large, there is no disk $D$ in $N(B_k^-)$ transverse to the $I$--fibers and with $\partial D\subset A_k$.  Now, suppose $D\subset N(B_k)$ is a disk in $S_n$ with $\partial D\subset A_k$.  So $D$ cannot be totally in $N(B_k^-)$.  If $\mu\cap D\ne\emptyset$ under any $B_k$--isotopy, since $\mu$ is the inverse limit of the infinite sequence of splittings, these splittings $\{B_k^-\}$ will eventually cut through $D$. By the claim above, there are only finitely many such disks $D$.  So, if $m$ is sufficiently large, there is no such disk $D\subset N(B_m)$ with $\mu\cap D\ne\emptyset$.  If $D\cap\mu=\emptyset$, since $D$ cannot be totally in $N(B_k^-)$ as above, we can split $B_k$ and $B_k^-$ further so that $D$ is carried by $B_k-B_k^-$ and hence $\partial D\not\subset A_k$ after this splitting.  Since there are only finitely many such disks $D$, after a finite sequence of splittings, we get a branched surface $B_k$ satisfying the requirements of the lemma.  

We should note that the assumption that $B$ does not carry any normal torus is important.  For example, if $\mu$ is a torus, one can easily construct a counter-example using an infinite sequence of disks wrapping around $\mu$ like the Reeb component.
\end{proof}

\begin{lemma}\label{Limmerse}
Let $M$, $B$, $\{S_n\}$ and $\mu$ be as in Lemma~\ref{Lessential}.  Let $\gamma$ be an immersed essential closed curve in a leaf of $\mu$.  Then $B$ can be split into a branched surface $B_1$ that carries $\mu$ and an infinite sub-sequence $\{S_{n_k}\}$ of $\{S_n\}$, such that, for each $k$, $S_{n_k}$ contains no embedded disk $D$ with the property that $\pi(\partial D)=\pi(\gamma)$, where $\pi: N(B_1)\to B_1$ is the collapsing map.
\end{lemma}
\begin{proof}
This lemma is basically the same as Lemma~\ref{Lessential}.  Although the curve $\gamma$ may not be embedded, each $S_n$ is embedded.  Hence there are only finitely many different configurations for $\partial D$.  So the lemma follows from the same arguments in the proof of Lemma~\ref{Lessential}.
\end{proof}

\section{Helix-turn-helix bands}\label{Shelix}

\begin{figure}
\begin{center}
\includegraphics{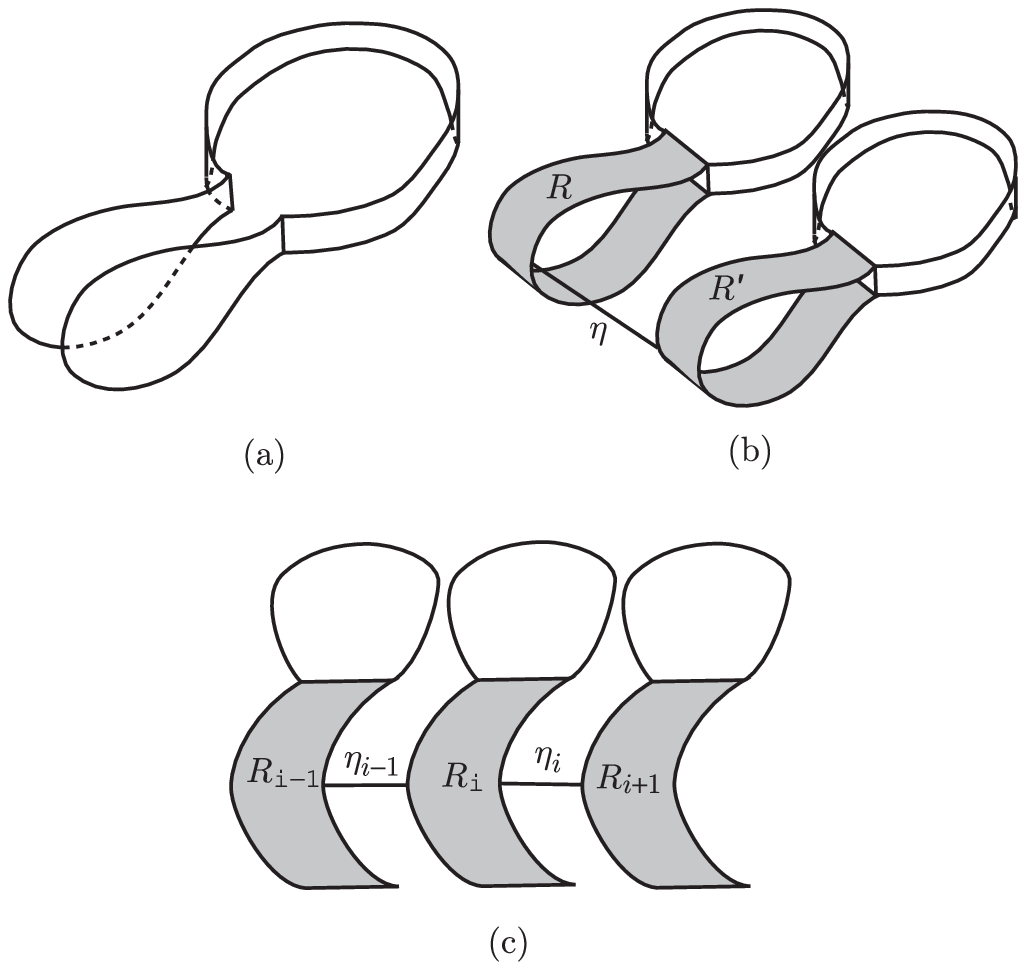}
\caption{}\label{monogon}
\end{center}
\end{figure}

A technical part in the proof of the main theorem is to construct compressing disks for the two handlebodies of the Heegaard splitting using $N(B)$.  Such compressing disks are constructed using a complicated band in $N(B)$ that connects two parallel monogons, as shown in Figure~\ref{monogon} (a).  The purpose of this section is to demonstrate how to construct these bands.  Such bands are constructed using a local picture of the limit lamination of a sequence of Heegaard surfaces.  We will start with a one-dimension lower example.

\begin{definition}
Let $A=S^1\times I$ be an annulus and $\alpha$ a compact spiral in $A$ transverse to the $I$--fibers.  We define the \emph{winding number} of $\alpha$, denoted by $w(\alpha)$, to be the smallest intersection number of $\alpha$ with an $I$--fiber of $A$.
\end{definition}

\begin{example}\label{exspiral}
Let $\tau$ be a train track obtained by attaching two ``tails" to a circle $\gamma$, as shown in Figure~\ref{traintrack}(b).  Curves fully carried by $\tau$ must consist of spirals. We use $x_1,\dots,x_4$ to denote the 4 branch sectors of $\tau$, and the branch equations are $x_1+x_3=x_4$ and $x_2+x_3=x_4$.  Suppose $\{\gamma_n\}$ is an infinite sequence of positive integer solutions to the branch equations whose limit $\mu$ in $\mathcal{PL}(\tau)$ is a measured lamination consisting of parallel circles carried by $\tau$.  So the coordinates of $\mu$ in $\mathcal{PL}(\tau)$ are $(0,0,1/2,1/2)$.  Let $\gamma_i=(x_1^{(i)},\dots,x_4^{(i)})\in\mathcal{S}(\tau)$ be the corresponding sequence of integer points.  For each $\gamma_i$, we denote the number of components of $\gamma_i$ by $h(\gamma_i)$ and clearly, $h(\gamma_i)=x_1^{(i)}=x_2^{(i)}$. Moreover, the winding number of each component of $\gamma_i$ is $w(\gamma_i)=x_3^{(i)}/h(\gamma_i)$.  Because of the branch equations, we have $\gamma_i=(x_1^{(i)},x_1^{(i)},x_3^{(i)}, x_1^{(i)}+x_3^{(i)})$.  Since the limit of these points in $\mathcal{PL}(\tau)$ is $(0,0,1/2,1/2)$, by part 1 of Proposition~\ref{Plinear}, we have that $\lim_{i\to\infty}x_1^{(i)}/x_3^{(i)}=0$, in other words $\lim_{i\to\infty}w(\gamma_i)=\infty$.  

In general, a train track near a circle can have many ``tails" like Figure~\ref{traintrack}~(a), but the argument above still works (using part 2 of Proposition~\ref{Plinear}).  If the limit of a sequence of spiral curves $\{\gamma_i\}$ is a measured lamination by circles, then the winding numbers tend to infinity, $\lim_{i\to\infty}w(\gamma_i)=\infty$.
\end{example}

Let $S^1\times I$ be an annulus, and let $\gamma$ be a collection of disjoint spirals properly embedded in $S^1\times I$ and transverse to the $I$--fibers. Suppose the winding number for each spiral is at least 2.  We fix an $I$--fiber $\{x\}\times I$.  Let $\beta$ be a subarc of a spiral in $\gamma$ with $\beta\cap(\{x\}\times I)=\partial\beta$. Let $\alpha$ be the subarc of $\{x\}\times I$ between the two endpoints of $\beta$.  We define the \emph{discrepancy} of $\gamma$ to be $1+|\gamma\cap int(\alpha)|$.  It is very easy to see that the discrepancy is equal to the number of components of $\gamma$ and does not depend on the choice of $\beta$.

Next, we consider the two-dimensional version of Example~\ref{exspiral}.

\begin{figure}
\begin{center}
\includegraphics{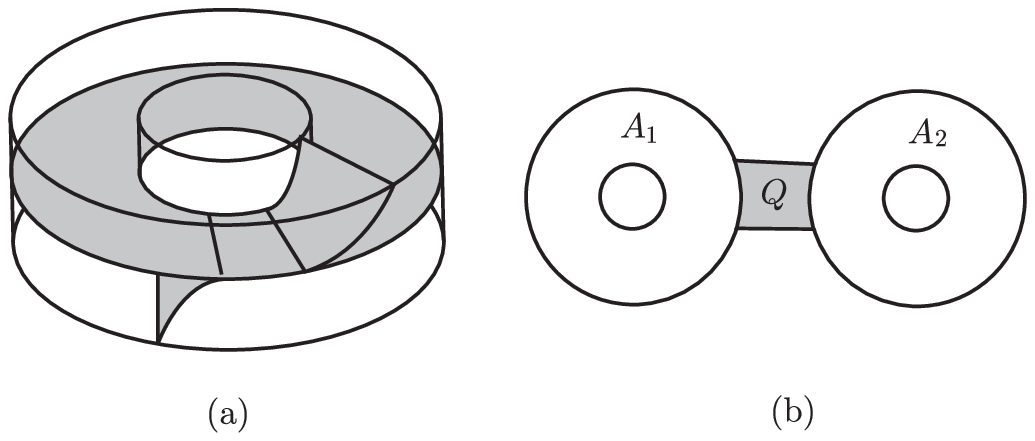}
\caption{}\label{cylinder}
\end{center}
\end{figure}

\begin{figure}
\begin{center}
\includegraphics{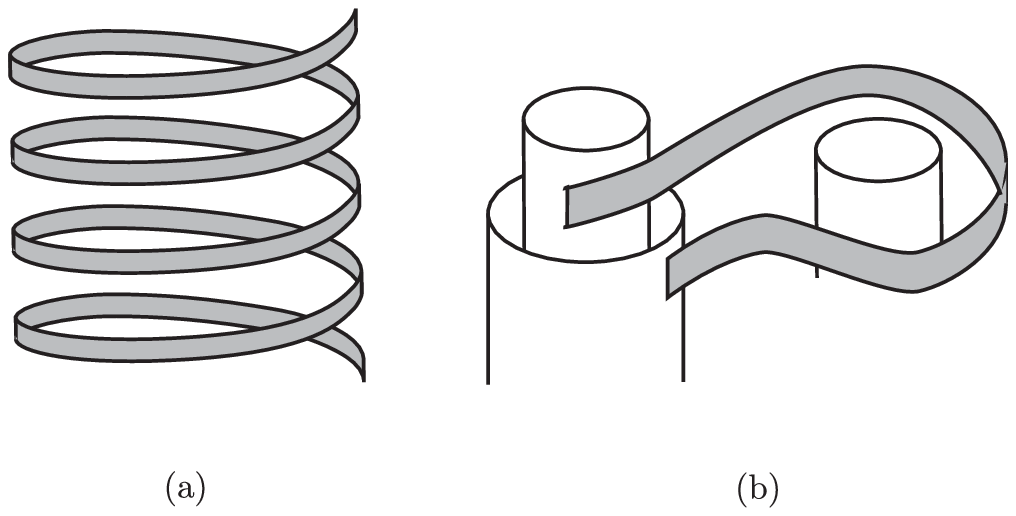}
\caption{}\label{band}
\end{center}
\end{figure}

\begin{example}\label{exband}
If we take a product of the train track in Example~\ref{exspiral} and an interval, we get a branched surface, see the shaded region in Figure~\ref{cylinder}~(a).  As in Figure~\ref{cylinder} (a), we may assume the branched surface is sitting in $A\times I$, where $A$ is a horizontal annulus, and this branched surface is transverse to the $I$--fibers of $A\times I$.  For any essential simple closed curve $c$ in $A$, the intersection of the cylinder $c\times I\subset A\times I$ and this branched surface is a train track as in Example~\ref{exspiral}.  Suppose there is a sequence of spiraling disks $\{S_n\}$ fully carried by this branched surface and the limit lamination of $\{S_n\}$ is a union of horizontal annuli of the form $A\times\{x\}$, $x\in I$.  Then we can define the winding number similarly, and if $n$ tends to infinity, the winding number of $S_n$ tends to infinity as well.  To fit this in the bigger picture, we should consider the $A\times I$ as a small portion of $N(B)$ and each $S_n$ is the intersection of $A\times I$ with a Heegaard surface.  Naturally, $S_n$ may not be connected.  Next we assume each $S_n$ lies in $A\times I$, transverse to every $I$--fiber of $A\times I$.  

Let $h$ be the number of components of $S_n$ and suppose $h\ge 2$.  Let $c$ be an essential simple closed curve in $A$.  We consider the vertical cylinder $c\times I\subset A\times I$.  $S_n\cap(c\times I)$ consists of $h$ spirals in $c\times I$. These spirals $S_n\cap(c\times I)$ cut $c\times I$ into some bands.  We may describe each band as a product $l\times J$,  where $l$ is a curve, $J$ is an interval, $l\times\partial J$ is a pair of spirals in $S_n\cap(c\times I)$, and each $\{x\}\times J$ ($x\in l$) is a subarc of an $I$--fiber of $c\times I$.  We call such a band $l\times J$ a \emph{helical band}, see the shaded region in Figure~\ref{band} (a) for a picture.  We call $\partial l\times J$ the two \emph{ends} of the band and define the wrapping number of the band to be the wrapping number of a spiral $l\times\{p\}$.  We define the \emph{thickness} of a helical band $l\times J$ to be the number of components of $S_n\cap(l\times J)$.  By the construction, the thickness of a helical band is at least 2 (since $l\times\partial J\subset S_n$) and can be as large as $h$.  If the thickness of a helical band is less than $h$, then we can find a larger helical band $l'\times J'$ that contains $l\times J$ and with larger thickness.  We say $l'\times J'$ is obtained by thickening $l\times J$.
\end{example}

\begin{example}\label{exjoint}
Let $A_1$ and $A_2$ be two annuli and $Q$ a quadrilateral.  By connecting $A_1$ and $A_2$ using $Q$, we get a pair of pants $P$, as shown in Figure~\ref{cylinder} (b).  Now we consider the product $P\times I$.  Let $\{S_n\}$ be a sequence of compact surface in $P\times I$ transverse to the $I$--fibers, and suppose the limit lamination of $\{S_n\}$ is of the form $P\times C$, where $C$ is a closed set in $I$.  Suppose each component of $S_n\cap(Q\times I)$ is of the form $Q\times\{x\}$, $x\in int(I)$, and suppose $S_n\cap (A_i\times I)$ ($i=1,2$) consists of spiraling disks as in Example~\ref{exband}.   We will use $h_i$ to denote the number of components of $S_n\cap(A_i\times I)$, and use $w_i$ to denote the winding number of a spiraling disk in $S_n\cap(A_i\times I)$.  As in Example~\ref{exspiral}, if $n$ is large, the winding number $w_i$ is large.  In this paper, we will also assume each $h_i$ is an even number and $h_i\ge 2$,  $i=1,2$.
\end{example}

$S_n\cap(A_i\times I)$ consists of $h_i$ spiraling disks ($i=1,2$).  Topologically, each spiraling disk is a meridian disk of the solid torus $A_i\times I$, and the intersection of theses spiraling disks with each annulus in $A_i\times\partial I$ is a union of parallel essential arcs in the annulus.  We say an arc $K$ is a \emph{proper vertical arc} if $K$ is a subarc of an $I$--fiber of $P\times I$ and $K$ is properly embedded in $\overline{P\times I-S_n}$.  Let $\gamma\times J$ be an embedded rectangle in $P\times I$. We call $\gamma\times J$ a \emph{vertical band} if each $\{x\}\times J$ is a subarc of an $I$--fiber and $\gamma\times\partial J$ lies in $S_n$.  Note that the helical bands described in Example~\ref{exband} are vertical bands.  We define the \emph{thickness} of the vertical band $\gamma\times J$ to be the number of components of $S_n\cap(\gamma\times I)$.  So the thickness of a vertical band is at least 2.

By our assumptions, the number of components of $S_n\cap(Q\times I)$ is roughly $w_1h_1=w_2h_2$.  Let $J$ be a proper vertical arc in $Q\times I$ and $\alpha_J$ be an arc in $Q\times I$ connecting a point in $int(J)$ to $Q\times\{0\}$.  We define the \emph{height} of $J$ to be the minimal number of intersection points in $S_n\cap\alpha_J$.  We take a vertical band $\beta\times J$ around $A_i\times I$ and with both vertical arcs $\partial\beta\times J$ in $Q\times I$, as in Figure~\ref{band} (b), and suppose each $\{x\}\times J$ ($x\in\beta$) is a proper vertical arc.  Then since $S_n\cap(A_i\times I)$ consists of spiraling disks, the height difference between the two proper vertical arcs $\partial\beta\times J$ is equal to the discrepancy (see the definition of discrepancy before Example~\ref{exband}) of the spirals around $A_i\times I$.  Hence the height difference between the two arcs $\partial\beta\times J$ is equal to $h_i$.  Moreover, two proper vertical arcs in $(A_i\cap Q)\times I$ belong to the same component of $\overline{A_i\times I-S_n}$ if and only if the height difference between the two arcs is $kh_i$ for some integer $k$.

Now we are in position to construct a helix-turn-helix band (the word helix-turn-helix comes from biology).  Recall that, as in Example~\ref{exjoint}, we assume each $h_i$ is an even number and $h_i\ge 2$.

\begin{example}[Helix-turn-helix bands]\label{HTH}
We assume $h_1=h_2$.  First we give an outline of the construction.  Let $c$ and $c'$ be a pair of disjoint essential simple closed essential curves in $A_1$. So $c\times I$ and $c'\times I$ are a pair of disjoint vertical annuli in $A_1\times I$.  We take a pair of helical bands in $c\times I$ and $c'\times I$ respectively and connect them using a vertical band going around $A_2\times I$, as depicted in Figure~\ref{band} (b).  The resulting vertical band is a helix-turn-helix band.  There are some subtleties and additional requirements. The detailed description of the construction is as follows.

Let $J_1$ be a proper vertical arc in $\overline{(A_1\cap Q)\times I}$.  We first take a vertical band $\sigma$ around $A_2\times I$, connecting $J_1$ to another proper vertical arc $J_2\subset (A_1\cap Q)\times I$, see the shaded region in Figure~\ref{band} (b).  Note that in Figure~\ref{band} (b), the left two cylinders are vertical cylinders in $A_1\times I$ and the right cylinder is a vertical cylinder in $A_2\times I$.  Clearly the height difference between $J_1$ and $J_2$ is $h_2$.  Since $h_1=h_2$, $J_1$ and $J_2$ lie in the same component of $\overline{A_1\times I-S_n}$.  Then we take a helical band $\sigma_i$ ($i=1,2$), as in Figure~\ref{band} (a), connecting $J_i$ to a proper vertical arc $J_i'$, where $J_i'$ has an endpoint in the bottom annulus $A_1\times\{0\}$.  We can choose $\sigma_1$ and $\sigma_2$ in different vertical cylinders in $A_1\times I$, see the left part of Figure~\ref{band} (b) for a picture of two disjoint cylinders.  So we may assume $\sigma_1\cap\sigma_2=\emptyset$ and $\Sigma=\sigma_1\cup\sigma\cup\sigma_2$ is an embedded vertical band connecting $J_1'$ to $J_2'$.  Note that since the height difference between $J_1$ and $J_2$ is $h_2=h_1$, the winding numbers for $\sigma_1$ and $\sigma_2$ differ by one.  We may write $\Sigma=\gamma\times J$, where $\gamma$ is an arc and $J$ is a closed interval.  $\Sigma$ has the properties that $\Sigma\cap S_n=\gamma\times\partial J$, each $\{x\}\times J$ is a subarc of an $I$--fiber of $P\times I$, and $\partial\gamma\times J=J_1'\cup J_2'$.  We call $\Sigma$ a \emph{helix-turn-helix} (or an HTH) band.  Note that the thickness of the vertical band $\Sigma$ in the construction above is 2.  Similar to Example~\ref{exband}, we can trivially thicken the HTH band $\Sigma$ to an embedded vertical band $\hat{\Sigma}$ so that the thickness of $\hat{\Sigma}$ is $h_1$ ($h_1=h_2$).  We call both $\Sigma$ and $\hat{\Sigma}$ HTH bands.

Since $J_1'$ and $J_2'$ lie in the same component of $\overline{A_1\times I-S_n}$ and each $J_i'$ has an endpoint in the bottom annulus $A_1\times\{0\}$, we may glue a small vertical band $\delta$ to $\Sigma$, connecting $J_1'$ to $J_2'$, and get a vertical annulus $A_\Sigma=\Sigma\cup\delta$ properly embedded in $\overline{(P\times I)-S_n}$.  Note that $J_1'\cup J_2'$ is a pair of opposite edges of $\delta$ and $\delta$ has an edge totally in the bottom annulus $A\times\{0\}$.  Let $x_i$ be the element in the fundamental group $\pi_1(P\times I)$ represented by the core of $A_i\times I$ ($i=1,2$).  Then this vertical annulus $A_\Sigma$ represents the element $x_1^{-k}\cdot x_2\cdot x_1^{k+1}$ in $\pi_1(P\times I)$, for some $k$.  
\end{example}

Note that in a previous version of the paper, there is a construction of an HTH band for the case $h_1<h_2$.  That construction turns out to be unnecessary for the proof of the main theorem.

\begin{example}\label{exmulti}
In Example~\ref{HTH}, if the winding numbers $w_1$ and $w_2$ are large, we can construct many disjoint HTH bands.  To see this, we first divide $P\times I$ into $N$ parts, $P\times I_i$ ($i=1,\dots,N$), where $I_i=[\frac{i-1}{N},\frac{i}{N}]$.  We may assume the intersection of $S_n$ with each $P\times I_i$ is as described in Example~\ref{exjoint}.  Suppose the winding numbers $w_1$ and $w_2$ are large. We can carry out the construction in Example~\ref{HTH} on each $P\times I_i$.  Then we glue a pair of long helical bands to the two ends of each vertical band constructed in $P\times I_i$ to spiral down to the bottom annulus $A_1\times\{0\}$.  By choosing these helical bands to be in disjoint vertical cylinders of $A_1\times I$ (see the left part of Figure~\ref{band} (b) for a picture of two disjoint cylinders), we may assume these HTH bands are disjoint.  Let $\Sigma_i=\gamma_i\times J$ ($i=1,\dots, N$) be the $N$ disjoint HTH bands above.  We may assume each component of $\partial\gamma_i\times J$ is a proper vertical arc with an endpoint in $A_1\times\{0\}$.  We may also construct the HTH bands so that these $\Sigma_i$'s lie in the same component of $\overline{P\times I-S_n}$. Moreover, we may assume that, for each $i$, the two proper vertical arcs $\partial\gamma_i\times J$ are close to each other.  Hence, similar to Example \ref{HTH}, we can glue a small vertical band $\delta_i$ to each $\Sigma_i$ and get a collection of disjoint vertical annuli $A_{\Sigma_i}=\Sigma_i\cup\delta_i$ ($i=1,\dots,N$) properly embedded in the same component of $\overline{P\times I-S_n}$.  The elements represented by these $A_{\Sigma_i}$'s in $\pi_1(P\times I)$ are conjugate. In fact, by unwinding the pairs of helical bands, we can isotope these annuli $A_{\Sigma_i}$ in $\overline{P\times I-S_n}$ so that $\pi(A_{\Sigma_i})$ is the same closed curve in $P$ for all $i$, where $\pi:P\times I\to P$ is the projection.  Furthermore, similar to Example~\ref{HTH}, we can trivially thicken these $\Sigma_i$'s into a collection of embedded disjoint HTH bands with thickness $h_1$.
\end{example}

Let $\Sigma$ and $A_\Sigma$ be the HTH band and the vertical annulus constructed in the examples above.  So, after a small perturbation, we may assume $\pi(A_\Sigma)$ is an immersed essential closed curve in $P$, where $\pi:P\times I\to P$ is the projection.  By Example~\ref{exmulti}, if $w_1$ and $w_2$ are large, we can choose $N$ disjoint HTH bands $\Sigma_i$ ($i=1,\dots, N$) and $N$ disjoint vertical annuli $A_{\Sigma_i}$.  Moreover, after some isotopy, $\pi(A_{\Sigma_i})$ is the same curve in $P$ for all $i$.  Thus, regardless of the configurations of $S_n$, as long as $n$ is large, there is a fixed finite set of immersed essential closed curves in $P$, denoted by $\mathcal{C}_P$, such that $\pi(A_{\Sigma_i})$ above is a curve in $\mathcal{C}_P$, up to isotopy.

The following lemma follows trivially from Lemma~\ref{Limmerse}.

\begin{lemma}\label{Lcor}
Let $M$, $B$, $\{S_n\}$ and $\mu$ be as in Lemma~\ref{Lessential}.  Let $P$ be an essential sub-surface of a leaf $l$ of $\mu$.  Suppose $P$ is a pair of pants.   Let $\mathcal{C}_P$ be the finite set of curves in $P$ as above.  Then $B$ can be split into a branched surface $B_1$ that carries $\mu$ and an infinite sub-sequence $\{S_{n_k}\}$ of $\{S_n\}$, such that no $S_{n_k}$ contains any disk $D$ with the property that $\pi(\partial D)=\pi(\gamma)$ for any $\gamma\in\mathcal{C}_P$.
\end{lemma}
\begin{proof}
By the hypotheses, every curve $\gamma\in\mathcal{C}_P$ is essential in the leaf $l$. So the lemma follows from Lemma~\ref{Limmerse}.
\end{proof}

\begin{definition}\label{Dmonogon}
Let $S_n$ be a closed embedded surface carried by $N(B)$, and let $\nu$ be a subarc of an $I$--fiber of $N(B)$ with $\partial\nu\subset S_n$.  We say that $\nu$ bounds a monogon if there is an embedded disk $E$ transverse to $S_n$, such that $\partial E=\nu\cup\alpha$, where $\alpha\subset S_n$ and $\partial\alpha=\partial\nu$.  We call the disk $E$ a \emph{monogon}, see Figure~\ref{yinyang}(b) for a picture.  We call $E$ an innermost monogon if $E\cap S_n=\alpha$.  Since $\nu\subset N(B)$, we may assume that a neighborhood of $\nu$ in $E$ is a sub-disk $\kappa=a\times J$ of $E$ such that each $\{x\}\times J$ is a subarc of an $I$--fiber of $N(B)$, $a\times\partial J\subset\alpha\subset\partial E$, and $\nu$ is a component of $\partial a\times J$.  We call $\kappa$ the \emph{tail} of the monogon.  We define the \emph{thickness} of the tail to be $|S_n\cap\nu|$ and define the length of the tail to be the length of a component of $a\times\partial J$.  So, if $E$ is innermost, the thickness of the tail is 2.   Let $\Sigma=\gamma\times J$ be an HTH band constructed in Example \ref{HTH}, and let $\nu$ be a component of $\partial\gamma\times J$.  Suppose $\nu$ bounds a monogon $E$ disjoint from $\Sigma$.  Then we can glue $\Sigma$ and two parallel copies of $E$ together, forming an embedded disk $\Delta$ as shown in Figure~\ref{monogon} (a).  By our construction, $\partial\Delta$ is a simple closed curve in $S_n$.  We call the disk $\Delta$ constructed in this fashion a \emph{pinched disk}.  Since $\Delta$ is constructed using parallel copies of $E$, there is a rectangle $R\subset S_n$ between the two monogons, see the shaded regions of Figure~\ref{monogon} (b).  Let $\Sigma'$ be another HTH band constructed in Example~\ref{exmulti}.  We can glue $\Sigma'$ and another two parallel copies of $E$ together, forming an embedded disk $\Delta'$.  Similarly, there is a rectangle $R'\subset S_n$ between the two monogons, as shown in Figure~\ref{monogon} (b).  By our construction in Example~\ref{exmulti}, $\Delta\cap\Delta'=\emptyset$ and $R\cap R'=\emptyset$.  Moreover, there is a short arc $\eta\subset S_n$ connecting $R$ to $R'$, as shown in Figure~\ref{monogon} (b).  We call an arc $\eta$ constructed in this fashion an \emph{$\eta$--arc}.
\end{definition}
\begin{remark}\label{Rrec}
Let $\Sigma$, $\Delta$ and $R$ be as in Definition~\ref{Dmonogon}.  We can denote $R=\alpha\times\beta$, where $\alpha$ and $\beta$ are intervals, and suppose $R\cap\Delta=\alpha\times\partial\beta\subset\partial\Delta$.  Moreover, $(\partial\Delta-\alpha\times\partial\beta)\cup(\partial\alpha\times\beta)$ is exactly the boundary of the annulus $A_\Sigma$ constructed in Example~\ref{HTH}.
\end{remark}

\begin{lemma}\label{Ldouble}
Let $M$, $\mu$, $P$, $B_1$ and $\{S_{n_k}\}$ be as in Lemma~\ref{Lcor}.  Suppose $P\times I$ is embedded in $N(B_1)$ with each $\{x\}\times I$ a subarc of an $I$--fiber of $N(B_1)$.  Suppose $S_{n_k}\cap(P\times I)$ is a surface as described in Example~\ref{exjoint} and assume the two winding numbers $w_1$ and $w_2$ are large enough.  Let $\Sigma=\gamma\times J$ be an HTH band constructed in the examples above.   Suppose the arcs $\partial\gamma\times J$ bound a pair of parallel embedded monogons $E_1$ and $E_2$ in $M-P\times (\epsilon,1]$, where $\epsilon\in I$ is a small number such that $\partial\gamma\times J\subset P\times[0,\epsilon]$.  As in Definition~\ref{Dmonogon},  let $\Delta=E_1\cup\Sigma\cup E_2$ be an embedded pinched disk with $\partial\Delta\subset S_{n_k}$.  Then $\partial\Delta$ is essential in $S_{n_k}$.
\end{lemma}
\begin{proof}
As in Example~\ref{HTH}, we can glue a small rectangle $\delta$ to $\Sigma$ and form an embedded annulus $A_\Sigma$ with $\partial A_\Sigma\subset S_{n_k}$.  By Lemma~\ref{Lcor}, $\partial A_\Sigma$ is a pair of essential curves in $S_{n_k}$.  Note that, since $E_1$ and $E_2$ may not be innermost monogons, $\Delta\cap S_{n_k}$ may contain other circles. 

Since $E_1$ and $E_2$ are parallel monogons, there is a thin rectangle $R\subset S_{n_k}$ between $\partial E_1$ and $\partial E_2$, and $E_1\cup E_2\cup\delta\cup R$ is an embedded 2--sphere in $M$, see the shaded region in Figure~\ref{monogon} (b) for a picture of $R$.

Since the two winding numbers $w_1$ and $w_2$ in Example~\ref{exjoint} are large, we may assume the number $\epsilon$ in the lemma is very small.  Hence, as in Example~\ref{exmulti}, we can find another disjoint HTH band $\Sigma'=\gamma'\times J$ and construct an annulus $A_{\Sigma'}$ by gluing a small rectangle $\delta'$ to $\Sigma'$.   By Lemma~\ref{Lcor}, $\partial A_{\Sigma'}$ is also essential in $S_{n_k}$.  Moreover, we can choose $\Sigma'$ so that $\partial\gamma'\times J$ bounds a pair of monogons $E_1'$ and $E_2'$ that are parallel to $E_1$ and $E_2$.  Similarly, $\Delta'=E_1'\cup\Sigma'\cup E_2'$ is also an embedded disk with $\partial\Delta'\subset S_{n_k}$ and $\Delta\cap\Delta'=\emptyset$.

Similar to $R$, there is also a thin rectangle $R'\subset S_{n_k}$ between $\partial E_1'$ and $\partial E_2'$.  Moreover, by our construction in Example~\ref{exmulti}, $R\cap R'=\emptyset$.  Since the 4 monogons $E_1$, $E_2$, $E_1'$, $E_2'$ are parallel to each other, as described in Definition~\ref{Dmonogon}, there is a short arc $\eta\subset S_{n_k}$ outside $P\times I$ connecting $R$ to $R'$, as shown in Figure~\ref{monogon} (b), where the two shaded regions are $R$ and $R'$.

Now, suppose $\partial\Delta$ is a trivial curve in $S_{n_k}$, and we use $D$ to denote the disk in $S_{n_k}$ bounded by $\partial\Delta$.  By Remark~\ref{Rrec}, $R\cap\partial\Delta$ is a pair of opposite edges of $R$, and the union of the other pair of opposite edges of $R$ and $\partial\Delta-R$ is $\partial A_\Sigma$.    Since $\partial A_{\Sigma}$ is essential in $S_{n_k}$, the rectangle $R$ must lie in $S_{n_k}-int(D)$.   Hence the arc $\eta$ must lie in $D$.  Since $R\cup\partial\Delta$ is disjoint from $R'\cup\partial\Delta'$, $R'\cup\partial\Delta'$ must lie in $D$.  This implies $\partial A_{\Sigma'}$ lies in $D$ and hence is trivial in $S_{n_k}$, contradicting our assumptions.
\end{proof} 

Therefore, after some splittings and taking a sub-sequence of $\{S_n\}$, we have the following.  For each HTH band $\Sigma$, by Lemma~\ref{Lcor}, the boundary of the annulus $A_\Sigma$ constructed above is a pair of essential curves in $S_n$.  Moreover, if the two ends of $\Sigma$ bound a pair of parallel monogons, by Lemma~\ref{Ldouble}, the boundary of the pinched disk $\Delta$ constructed above is also an essential curve in $S_n$. 

\section{Proof of the main theorem}\label{Smain}

Suppose $M$ is a closed orientable irreducible and non-Haken 3--manifold and $M$ is not a Seifert fiber space.  By Theorem~\ref{Heeg1}, $M$ has a finite collection of branched surfaces such that, 
\begin{enumerate}
\item each branched surface in this collection is obtained by gluing together normal disks and at most one almost normal piece with respect to a fixed triangulation, similar to \cite{FO},
\item up to isotopy, every strongly irreducible Heegaard surface is fully carried by a branched surface in this collection. 
\item no branched surface in this collection carries any normal 2--sphere or normal torus.
\end{enumerate}
The goal of this section is to prove Theorem~\ref{Tfinite}.  It is clear that Theorem~\ref{Tfinite} and Theorem~\ref{Heeg1} imply the main theorem.

\begin{theorem}\label{Tfinite}
Suppose $M$ is a closed, orientable, irreducible and non-Haken 3--manifold.  Let $B$ be a branched surface in Theorem~\ref{Heeg1}.  Then $B$ carries only finitely many irreducible Heegaard surfaces, up to isotopy.
\end{theorem}
\begin{proof}
Each closed surface fully carried by $B$ corresponds to a positive integer solution to the branch equations.  Since the projective lamination space $\mathcal{PL}(B)$ is compact, if $B$ fully carries an infinite number of distinct strongly irreducible Heegaard surfaces, then there is an accumulation point in the projective lamination space, which corresponds to a measured lamination $\mu$.  We may consider $\mu$ as the limit of these Heegaard surfaces, see section~\ref{Slimit}.  Our goal is to show that $\mu$ is also an essential lamination. Then by Proposition~\ref{PHaken}, $M$ is Haken, which contradicts our hypothesis.

Because of Theorem~\ref{TMS}, we divide the proof into two parts.  Part A is the case that $\mu$ is an exceptional minimal lamination and part B is the case that $\mu$ is a closed surface. The proofs for the two cases are slightly different.  

\medskip

\noindent\underline{\textbf{Part A}}. $\mu$ is an exceptional minimal measured lamination.

\medskip

The main task is to prove the following lemma.

\begin{lemma}\label{Lincomp}
$\mu$ is incompressible in $M$. 
\end{lemma}
\begin{proof}[Proof of Lemma~\ref{Lincomp}]
Suppose $\{S_n\}$ is an infinite sequence of strongly irreducible Heegaard surfaces fully carried by $B$ and $\mu$ is the limit point of $\{S_n\}$ in $\mathcal{PL}(B)$.  The lamination $\mu$ is carried by $B$, but it may not be fully carried by $B$.  Let $B^-$ be the sub-branched surface of $B$ that fully carries $\mu$. By Corollary~\ref{Cnormal}, $\mu$ must be a normal lamination.  Hence $B^-$ does not contain the almost normal piece and every surface carried by $B^-$ is normal.  By our hypotheses, $B^-$ does not carry any 2--sphere or torus.  

We may assume $N(B^-)\subset N(B)$ with the induced $I$--fiber structure.  By Proposition~\ref{Psplit}, we can arbitrarily split $B^-$ along $\mu$ and then split $B$ accordingly so that the resulting branched surface still carries an infinite sub-sequence of $\{S_n\}$.  Therefore, by Proposition~\ref{Psplit} and Lemma~\ref{Lfine}, after splitting $B$ and $B^-$ and taking an infinite sub-sequence of $\{S_n\}$, we may assume no component of $\partial_hN(B^-)$ is a disk and each leaf of $\mu$ is $\pi_1$--injective in $N(B^-)$.

After some isotopy, we may assume $\partial_hN(B^-)\subset\mu$.  Suppose $\mu$ is compressible and let $D$ be a compressing disk.  After some splittings on $B$ and $B^-$ as in Proposition~\ref{Psplit} and taking a sub-sequence of $\{S_n\}$, we may assume $\partial_hN(B^-)$ is compressible and $D$ is a compressing disk in $M-int(N(B^-))$.

So $\gamma_1=\partial D$ is an essential curve in a leaf $l$ of $\mu$.  Since $\mu$ has no holonomy, there is a vertical annulus $V$ in $N(B^-)$ such that $V$ contains $\gamma_1$ and $\mu\cap V$ is a union of parallel circles.  Thus, after some splittings on $B^-$, we may assume $\pi(\gamma_1)$ is a simple closed curve in $B^--L^-$ and $V=\pi^{-1}(\pi(\gamma_1))$, where $\pi: N(B^-)\to B^-$ is the collapsing map and $L^-$ is the branch locus of $B^-$.   By Proposition~\ref{Psplit} and Remark~\ref{Rsss}, we may split $B$ accordingly and assume $B$ still carries an infinite sequence of Heegaard surfaces $\{S_n\}$ whose limit lamination is $\mu$ and $B^-$ is the sub-branched surface of $B$ that fully carries $\mu$.

By Lemma~\ref{Lessential}, after some splittings and taking a sub-sequence of $\{S_n\}$, we may assume that $S_n\cap V$ does not contain any circle that is trivial in $S_n$, for each $n$.  Since $\gamma_1$ bounds an embedded disk in $M$, by Lemma~\ref{Ls22}, if $S_n\cap V$ consists of circles, then each circle bounds a compressing disk in one of the two handlebodies.  However, if $S_n\cap V$ consists of circles, by Corollary~\ref{Clinear} and Example~\ref{exspiral},  the number of circles in $S_n\cap V$ tends to infinity as $n$ tends to infinity.  This gives a contradiction to Lemma~\ref{Lbound}.  Therefore, $S_n\cap V$ cannot be a union of circles if $n$ is large enough.  So we may assume $S_n\cap V$ consists of spirals.

Since every leaf is dense, $l\cap V$ contains an infinite number of circles.  Since $B^-$ does not carry any torus, by Proposition~\ref{Ptorus} and our assumptions on $N(B^-)$ above, there must be a circle $\gamma_2\subset l\cap V$ such that $\gamma_2$ is non-trivial and not homotopic to $\gamma_1$ in $l$.  Let $\gamma_i\times I\subset V$ ($i=1,2$) be a pair of disjoint thin vertical annuli such that $\gamma_i\subset\gamma_i\times I$ and $\mu\cap(\gamma_i\times I)$ is a union of parallel circles.  Let $\alpha\subset l$ be a simple arc connecting $\gamma_1$ to $\gamma_2$ and $\Gamma=\gamma_1\cup\alpha\cup\gamma_2$ be a 1--complex in $l$.  By choosing $\gamma_i\times I$ to be thin enough, we may assume $\gamma_1\times I$ and $\gamma_2\times I$ are connected by a rectangle $\alpha\times I$, forming an embedded 2--complex $\Gamma\times I$ with each $\{x\}\times I$ ($x\in\Gamma$) a subarc of an $I$--fiber of $N(B^-)$.  By our construction, $\mu\cap(\Gamma\times I)$ is a union of 1--complexes parallel to $\Gamma$.

Let $A_i\subset l$ ($i=1,2$) be a small annular neighborhood of $\gamma_i$ in $l$, $Q$ be a small neighborhood of $\alpha$ in $l$, and $P=A_1\cup Q\cup A_2$ be a small neighborhood of $\Gamma$ in $l$. We can extend $\Gamma\times I$ to a product $P\times I\subset N(B^-)$.  So $\mu\cap(P\times I)$ is a union of compact surfaces parallel to $P$.  Moreover, since $\gamma_1$ and $\gamma_2$ are not homotopic in $l$, $P$ is an essential sub-surface of $l$.   

Since every leaf is dense in $\mu$, after some splittings along $\mu$, we may assume $\pi(P)\subset B^--L^-$ and $P\times I=\pi^{-1}(\pi(P\times I))$, where $\pi: N(B^-)\to B^-$ is the collapsing map and $L^-$ is the branch locus of $B^-$.  By Proposition~\ref{Psplit}, we may split $B$ accordingly and assume $B$ still carries an infinite sequence of Heegaard surfaces $\{S_n\}$ whose limit lamination is $\mu$ and $B^-$ is the sub-branched surface of $B$ that fully carries $\mu$. 

By the construction above, we may consider $\mu\cap(P\times I)$ as the limit lamination of the sequence $\{S_n\cap(P\times I)\}$.  Since $S_n\cap V$ consists of spirals and $\gamma_i\times I\subset V$, after some splittings, we may assume $S_n\cap(A_i\times I)$ is a union of spiraling disks and $S_n\cap(P\times I)$ is as described in Example~\ref{exjoint}.  We use the same notations as section~\ref{Shelix}, and in particular, let $h_i$ be the number of components of $S_n\cap(A_i\times I)$.  Since $\gamma_1\times I$ and $\gamma_2\times I$ are disjoint sub-annuli of $V$ before the splitting, we may assume the spirals in $S_n\cap V$ wind around both $\gamma_i\times I$ many times. So the spirals in $S_n\cap(\gamma_i\times I)$ are part of longer spirals in $S_n\cap V$.  Hence, the discrepancies (see section~\ref{Shelix} for the definition of discrepancy) of the spirals in $S_n\cap(\gamma_i\times I)$ ($i=1,2$) are the same.  Therefore, we have $h_1=h_2$.

By Lemma~\ref{Lcor}, after some splittings and taking a sub-sequence of $\{S_n\}$, we may assume that no $S_n$ contains a disk $E$ with the property that $\pi(\partial E)=\pi(\gamma)$ for any curve $\gamma\in\mathcal{C}_P$, where $\mathcal{C}_P$ is as in Lemma~\ref{Lcor}.

Recall that $\gamma_1$ bounds a compressing disk $D$ in $M-int(N(B^-))$.  Let $\hat{D}=D\cup(\gamma_1\times I)$, where $\gamma_1\times I\subset P\times I$.  By our construction above, $\hat{D}$ is an embedded disk in $M$.   As $S_n$ is a compact surface, $S_n\cap \hat{D}$ must produce a monogon with a long ``tail" spiraling around $\gamma_1\times I$, as shown in Figure~\ref{yinyang}~(a).  In fact, there are at least 2 monogons as the $\theta_1$ and $\theta_2$ in Figure~\ref{yinyang}~(a).  Note that $\hat{D}\cap S_n$ may contain circles.

\begin{figure}
\begin{center}
\includegraphics{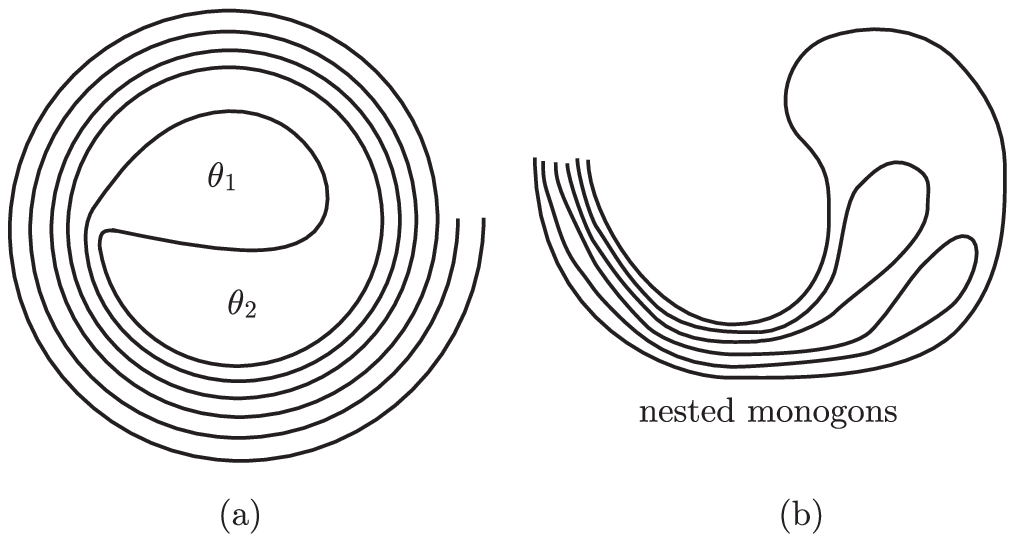}
\caption{}\label{yinyang}
\end{center}
\end{figure}

We will first consider the case that $\hat{D}\cap S_n$ does not contain any circle.  

\medskip

\noindent\underline{\textbf{Case 1}}. $\hat{D}\cap S_n$ does not contain any circle. 

\medskip

Since $S_n$ is a separating surface, we may assume each $h_i$ (i.e. the number of components in $S_n\cap(A_i\times I)$) is an even number.  Since $h_1=h_2$, we have the following 2 subcases.

\medskip

\noindent\underline{\textbf{Subcase 1a}}.  $h_1=2$.

\medskip

In this case, $S_n\cap\hat{D}$ is basically a single curve with both ends wrapping around $\gamma_1\times I$, as shown in Figure~\ref{yinyang}(a).  So we have two innermost monogons $\theta_1$ and $\theta_2$ in different handlebodies.  After a small perturbation in a small neighborhood of $\hat{D}$, we may assume $\partial\theta_1$ and $\partial\theta_2$ are disjoint in $S_n$.  As in Examples~\ref{exmulti}, we can find two disjoint HTH bands $\Sigma_1$ and $\Sigma_2$, such that $\Sigma_i$ ($i=1,2$) connects two parallel copies of $\theta_i$ forming a pinched disk $\Delta_i$.  By Lemma~\ref{Ldouble} and our construction, $\Delta_1$ and $\Delta_2$ are compressing disks in different handlebodies and $\partial\Delta_1\cap\partial\Delta_2=\emptyset$, which contradicts the assumption that $S_n$ is a strongly irreducible Heegaard surface.

\medskip

\noindent\underline{\textbf{Subcase 1b}}.  $h_1\ge 4$.

\medskip

Since $h_1\ge 4$, $\hat{D}\cap S_n$ contains at least two curves.  Note that each curve of $\hat{D}\cap S_n$ cuts $\hat{D}$ into 2 monogons, as the $\theta_1$ and $\theta_2$ in Figure~\ref{yinyang}~(a).  Thus, we can find a monogon $E$ which is not innermost, but each monogon in the interior of $E$ is innermost, as shown in Figure~\ref{yinyang}(b).  Let $h_E$ be the thickness of the tail of $E$ and clearly $h_E\le h_1$. 

Let $E_1$ and $E_2$ be two parallel copies of $E$.  Since $h_E\le h_1$, by Example~\ref{HTH}, we can connect the tails of $E_1$ and $E_2$ using an HTH band $\hat{\Sigma}$ (in $P\times I$) with thickness $h_E$.  We denote $E_1\cup\hat{\Sigma}\cup E_2$ by $\Delta$.  So $\Delta$ is an embedded disk with $\partial\Delta\subset S_n$.  Let $c_1,\dots,c_m$ be the components of $S_n\cap int(\Delta)$, and let $\Delta_i$ be the disk in $\Delta$ bounded by $c_i$.  Similar to $\Delta$, each $\Delta_i$ is the union of a sub-band of $\hat{\Sigma}$ and two parallel copies of a sub-monogon of $E$.  Since we have assumed each sub-monogon in $int(E)$ is innermost, $\Delta_i\cap S_n=\partial\Delta_i$.  By Lemma~\ref{Ldouble} and our assumptions above, $\partial\Delta$ and each $\partial\Delta_i$ are essential curves in $S_n$.  So each $\Delta_i$ is a compressing disk in a handlebody, say $H_1$.  By Lemma~\ref{Ls22}, $\partial\Delta$ must bound a disk in a handlebody.  Since the Heegaard splitting is strongly irreducible, $\partial\Delta$ must bound a compressing disk in $H_1$ as well.  So $P_\Delta=\Delta-\cup_{i=1}^mint(\Delta_i)$ is a planar surface properly embedded in $H_2$.  If $P_\Delta$ is compressible in $H_2$, then we can compress $P_\Delta$ into a collection of disjoint incompressible planar surfaces $P_1,\dots, P_s$.  By Corollary~\ref{Csch}, each $P_i$ is $\partial$--parallel in $H_2$.  Let $Q_i$ be the sub-surface of $S_n$ that is parallel to $P_i$ in $H_2$ ($\partial P_i=\partial Q_i$).  Since the $P_i$'s are disjoint, any two surfaces $Q_i$ and $Q_j$ are either disjoint or nested in $S_n$.  

By Example~\ref{exmulti}, we can construct another HTH band $\hat{\Sigma}'$ with thickness $h_E$, connecting two monogons $E_1'$ and $E_2'$, where $E_1'$ and $E_2'$ are also two parallel copies of $E$.  We use $\Delta'$ to denote the disk $E_1'\cup\hat{\Sigma}'\cup E_2'$.  By our construction $\Delta\cap\Delta'=\emptyset$.  Similar to $\Delta$, $int(\Delta')\cap S_n$ is a union of circles $c_1',\dots,c_m'$ and the sub-disk of $\Delta'$ bounded by $c_i'$, denoted by $\Delta_i'$, is a compressing disk for the handlebody $H_1$.  Similarly, we can compress the planar surface $\Delta'-\cup_{i=1}^mint(\Delta_i')$ into a collection of incompressible planar surfaces $P_1',\dots,P_t'$.  Since $\Delta\cap\Delta'=\emptyset$, we may assume these $P_i$'s and $P_j'$'s are all disjoint in $H_2$.  So each $P_i'$ is also $\partial$--parallel in $H_2$ and we use $Q_i'$ to denote the sub-surface of $S_n$ that is parallel to $P_i'$ ($\partial P_i'=\partial Q_i'$).  Since these planar surfaces $P_i$'s and $P_j'$'s are disjoint and $\partial$--parallel, any two surfaces $Q_i$ and $Q_j'$ are either disjoint or nested in $S_n$.  

To unify notations, we also denote $\Delta$, $\partial\Delta$, $\Delta'$, $\partial\Delta'$ by $\Delta_0$, $c_0$, $\Delta_0'$, $c_0'$ respectively.

As in Definition~\ref{Dmonogon} and Remark~\ref{Rrec}, for each $c_i=\partial\Delta_i$ (resp. $c_i'$), $0\le i\le m$, there is a rectangle $R_i=\alpha_i\times\beta_i$ (resp. $R_i'=\alpha_i'\times\beta_i'$) in $S_n$, see the shaded region in Figure~\ref{monogon}(b), such that $R_i\cap c_i$ (resp. $R_i'\cap c_i'$) is a pair of opposite edges $\alpha_i\times\partial\beta_i$ (resp. $\alpha_i\times\partial\beta_i$). Moreover, these $R_i$ and $R_j'$ are pairwise disjoint.  By our construction in section~\ref{Shelix}, $(c_i-\alpha_i\times\partial\beta_i)\cup(\partial\alpha_i\times\beta_i)$ (resp. $(c_i'-\alpha_i'\times\partial\beta_i')\cup(\partial\alpha_i'\times\beta_i')$) is the boundary of an embedded vertical annulus $A_{\Sigma_i}$ (resp. $A_{\Sigma_i'}$) in $P\times I$, and by our assumptions and Lemma~\ref{Lcor}, $\partial A_{\Sigma_i}$ (resp. $\partial A_{\Sigma_i'}$) is a pair of essential curves in $S_n$. 

Let $W_i$ (resp. $W_i'$) be the closure of a small neighborhood of $c_i\cup R_i$ (resp. $c_i'\cup R_i'$) in $S_n$.  So two boundary circles of $W_i$ (resp. $W_i'$) are parallel to the two components of $\partial A_{\Sigma_i}$ (resp. $\partial A_{\Sigma_i'}$) above, and the other boundary component of $W_i$ (resp. $W_i'$) is parallel to $c_i$ (resp. $c_i'$).  By our assumptions above, each boundary circle of $W_i$ (resp. $W_i'$) is an essential curve in $S_n$.  Moreover, there is an $\eta$--arc (see Definition~\ref{Dmonogon}) $\eta_i\subset S_n$ connecting $R_i$ to $R_i'$, as shown in Figure~\ref{monogon} (b).

Let $Q_i$ be a planar surface above, and suppose $c_0,\dots,c_q$ are the boundary components of $Q_i$.  Next, we will show that at least one $R_j$ ($0\le j\le q$) lies in $S_n-int(Q_i)$.  Otherwise, suppose $R_j\subset Q_i$ for every $j$.  Then for each $j$, $\partial\alpha_j\times\beta_j$ is a pair of arcs properly embedded in $Q_i$.  Since $Q_i$ is a planar surface and since there is a rectangle $R_j$ attached to each $c_j$, by an innermost-surface argument, it is easy to see that, for some $j$, $\partial\alpha_j\times\beta_j$ is a pair of $\partial$--parallel arcs in $Q_i$.  This implies that a boundary component of $W_j$ bounds a disk in $Q_i$ and hence is trivial in $S_n$, contradicting our assumptions.  This argument also holds for each $Q_i'$.  Therefore, for each $Q_i$ (resp. $Q_i'$), there is always such a rectangle $R_j$ (resp. $R_k'$), lying outside $int(Q_i)$ (resp. $int(Q_i')$) and with two opposite edges in $\partial Q_i$ (resp. $\partial Q_i'$).

Let $Q_i$ be any planar surface above. Suppose $c_k$ is a boundary circle of $Q_i$ and suppose $R_k=\alpha_k\times\beta_k$ is a rectangle outside $int(Q_i)$.  So $R_k\cap Q_i=R_k\cap c_k=\alpha_k\times\partial\beta_k$.  By our construction before, there is an arc $\eta_k$ connecting $R_k$ to $R_k'$, and $int(\eta_k)$ is disjoint from any $c_j$ or $c_j'$.  Moreover, the two endpoints of $\eta_k$ lie in $\alpha_k\times\partial\beta_k\subset  c_k$ and $\alpha_k'\times\partial\beta_k'\subset c_k'$.  Suppose $c_k'$ is a boundary component of $Q_j'$.  Since $R_k$ lies outside $int(Q_i)$, $\eta_k$ must lie in $Q_i$.  Hence $c_k'\subset Q_i$.  Since the planar surfaces $Q_i$ and $Q_j'$ are either disjoint or nested, $c_k'\subset Q_i$ implies that $Q_j'\subset Q_i$.  This means that for each $Q_i$, there is some $Q_j'$ such that $Q_j'\subset Q_i$.

However, we can apply the same argument to $Q_i'$ and conclude that, for each $Q_i'$, there is some $Q_k$ such that $Q_k\subset Q_i'$.  This is impossible because there is always an innermost planar surface among these $Q_i$'s and $Q_j'$'s.

\medskip

\noindent\underline{\textbf{Case 2}}. $\hat{D}\cap S_n$ contains circles.

\medskip

Similar to Case 1, each non-circular curve cuts $\hat{D}$ into a pair of monogons, though there may be circles in the monogons.  We say a monogon $E$ is innermost, if $E$ does not contain other monogon, but $E$ may contain circles of $\hat{D}\cap S_n$.  We first consider innermost monogons.  Let $E$ be an innermost monogon and $c_1,\dots, c_K$ the outermost circles of $E\cap S_n$.  Since the sequence of surfaces $\{S_n\}$ are carried by $B$, by assuming $\hat{D}$ to be transverse to $B$, it is easy to see that $K$, the number of such outermost circles in $E$, is bounded from above by a number independent of $S_n$.  Since we assume $n$ is large, the winding number $w_i$ of the spiraling disks in $A_i\times I$ is large.  So, by Example~\ref{exmulti}, we can find a large number of disjoint HTH bands $\Sigma_1,\dots,\Sigma_N$. Moreover, we can take $2N$ parallel copies of $E$, denoted by $E_1, E_1',\dots, E_N, E_N'$, so that the disks $\Delta_i=E_i\cup\Sigma_i\cup E_i'$ are disjoint and embedded in $M$.  By Lemma~\ref{Ldouble}, we may assume each $\partial\Delta_i$ is an essential curve in $S_n$.  Since $K$ is bounded by a number independent of $S_n$, we may assume $N$ is much larger than $K$, and this is a key point in the proof.

Between each pair $E_i$ and $E_i'$, there is a rectangle $R_i\subset S_n$ with two opposite edges in $\partial E_i$ and $\partial E_i'$, see the shaded region in Figure~\ref{monogon} (b).  By the construction in section~\ref{Shelix}, we may assume there is an $\eta$--arc (see Definition~\ref{Dmonogon}) $\eta_i$ connecting $R_i$ to $R_{i+1}$ for each $i=1,\dots,N-1$, as shown in Figure~\ref{monogon} (c).  The interior of each $\eta_i$ is disjoint from these disks $\Delta_j$'s.  

If $c_i$ ($i=1,\dots, K$) is a trivial curve in $S_n$, since $M$ is irreducible, we can perform some isotopy on $E$ (fixing $\partial E$) and get a monogon disk with fewer outermost circles in $E\cap S_n$.  So we may assume each $c_i$ is essential in $S_n$.  Let $d_i$ be the disk in the monogon $E$ bounded by $c_i$ ($i=1,\dots,K$), and suppose $E-\cup_{i=1}^Kd_i$ lies in $H_1$.  By Lemma~\ref{Ls22}, each circle $c_i$ bounds a compressing disk in a handlebody.  If some $c_i$  bounds a disk in $H_1$, then we can replace $d_i$ by a disk in $H_1$ and obtain a disk with the same boundary $\partial\Delta_i$ but fewer outermost circles.  If we can eliminate all the outermost circles $c_i$'s in this fashion, then we can conclude that each $\partial\Delta_i$ bounds a compressing disk in $H_1$.  Suppose we cannot eliminate these circles $c_i$ ($i=1,\dots,K$) via these isotopies and surgeries.  Then by the arguments above, we may assume each $c_i$ bounds a compressing disk in $H_2$. 

The arguments next involve compression bodies and strongly irreducible Heegaard splittings for manifold with boundary.  We refer to \cite{CG} for definitions and fundamental results. 

Let $W$ be the 3--manifold obtained by adding $K$ 2--handles to $H_1$ along these $c_i$'s, and let $\hat{W}$ be the manifold obtained by capping off the 2--sphere components of $\partial W$ by 3--balls.  Since each $\Delta_i$ is constructed using parallel copies of $E$, after some isotopies, we may assume each $\Delta_i$ is a properly embedded disk in $W$.  Note that after pushing $S_n$ into $int(\hat{W})$, $S_n$ becomes a Heegaard surface for $\hat{W}$, bounding the handlebody $H_1$ on one side and a compression body $W_2$ on the other side.  Since each $c_i$ bounds a compressing disk in $H_2$ and $M=H_1\cup_{S_n} H_2$ is strongly irreducible, by \cite{CG}, the Heegaard splitting $\hat{W}=H_1\cup_{S_n} W_2$ is also strongly reducible.  By our assumption on $E$ above, at least one 2--handle is added to $H_1$ and hence $W_2$ is not a trivial compression body.  Thus, by a theorem of Casson and Gordon (Theorem 2.1 of \cite{CG}), $\hat{W}$ is irreducible, and if $\partial\hat{W}\ne\emptyset$, $\partial\hat{W}$ is incompressible in $\hat{W}$.  Therefore, each $\partial\Delta_i$ bounds a disk $D_i$ in $\partial W$ ($i=1,\dots,N$).

Since $W$ is obtained by attaching $K$ 2--handles to $H_1$, there are $2K$ disjoint disks $m_1,\dots m_{2K}$ in $\partial W-S_n$ parallel to the cores of these 2--handles. Note that one can obtain the handlebody $H_2$ by attaching 1--handles to $M-int(W)$ along these disks $m_i$'s.  Since each $\partial\Delta_i=\partial D_i$ is essential in $S_n$, each disk $D_i\subset\partial W$ ($i=1,\dots,N$) must contain some $m_j$ ($1\le j\le 2K$).  Recall that $K$ is bounded by a number independent of $S_n$ and we have assumed that $N$ is very large compared with $K$.  Since each $D_i$ contains some $m_j$, for any integer $p$, if $N$ is large enough, there exist a sequence of $p$ nested disks $D_{a_1}\subset D_{a_2}\subset\cdots\subset D_{a_p}$ ($0\le a_i\le N$).  Note that if $p>2K$, at least one annulus $D_{a_{i+1}}-D_{a_i}$ does not contain any disk $m_j$.  So, by assuming $N$ is large enough, one can find 3 nested disks, say $D_1\subset D_2\subset D_3$, such that the two annuli $D_3-int(D_2)$ and $D_2-int(D_1)$ do not contain any disk $m_i$.  

Recall that $\partial D_2=\partial\Delta_2$ and there is a rectangle $R_2\subset S_n$ with two opposite edges attached to $\partial\Delta_2$.  By the construction of $W$, we also have $R_2\subset\partial W$.  Moreover, $int(R_2)$ is disjoint from the circles $\partial\Delta_j$'s.  So $R_2$ lies in one of the two annuli, $D_3-int(D_2)$ or $D_2-int(D_1)$.  Let $W(\partial\Delta_2\cup R_2)$ be the closure of a small neighborhood of $\partial\Delta_2\cup R_2$ in $S_n$. By our assumptions before, each boundary circle of $W(\partial\Delta_2\cup R_2)$ is essential in $S_n$.  Since the two annuli $D_3-int(D_2)$ and $D_2-int(D_1)$ do not contain any disk $m_i$, one boundary circle of $W(\partial\Delta_2\cup R_2)$ must be a trivial circle in both $\partial W$ and $S_n$, which contradicts our constructions and assumptions on the $R_i$'s before.

The arguments above show that, for any innermost monogon $E$ and pinched disk $\Delta_i$ above, after some isotopies and surgeries, we can eliminate these outermost circles $c_i$'s so that $\Delta_i$ becomes a compressing disk in $H_1$, where $E-\cup_{i=1}^Kd_i\subset H_1$ as above.  Now, similar to Case 1, we have two subcases.

\medskip

\noindent\underline{\textbf{Subcase 2a}}. $h_1=2$.

\medskip

In this subcase, $\hat{D}\cap S_n$ contains exactly one non-circular curve and this curve cuts $\hat{D}$ into a pair of innermost monogons.  So, by the arguments above on innermost monogons, we can eliminate the outermost circles in $\hat{D}\cap S_n$, and construct two disjoint compressing disks in the two handlebodies as in Subcase 1a.

\medskip

\noindent\underline{\textbf{Subcase 2b}}. $h_1\ge 4$.

\medskip

The proof for this subcase is a combination of the proof of Subcase 1b and the arguments on innermost monogons above.  Similar to Subcase 1b, we can find a monogon $E$ which is not innermost, but each monogon in the interior of $E$ is innermost.  As in Subcase 1b, by connecting two parallel copies of $E$ and a thick HTH band $\Sigma$, we get a pinched disk $\Delta$ (see Definition~\ref{Dmonogon}), with $\partial\Delta\subset S_n$.  Let $\epsilon_1,\dots,\epsilon_k$ be the monogons in $int(E)$.  Then the corresponding parallel copies of $\epsilon_i$ and a sub-band of $\Sigma$ form a pinched disk $\Delta_i\subset\Delta$ ($i=1,\dots,k$).  By Lemma~\ref{Ldouble} and our assumptions before, $\partial\Delta$ and each $\partial\Delta_i$ are essential in $S_n$.  By the arguments on innermost monogons, after some isotopies and surgeries, we may assume $S_n\cap int(\Delta_i)=\emptyset$ and each $\Delta_i$ is a compressing disk in a handlebody.  Since $S_n$ is strongly irreducible, these $\Delta_i$'s are compressing disks in the same handlebody, say $H_2$.

Let $c_1,\dots, c_{K}$ be the outermost circles in $E\cap S_n$.  As before, $K$ is bounded by a number independent of $S_n$. By our assumption on innermost monogons, these $c_i$'s lie in $E-\cup_{i=1}^k\epsilon_i$.  By the construction of the pinched disks, $S_n\cap(\Delta-\cup_{i=1}^k\Delta_i)$ has $2K$ outermost circles $c_1,\dots, c_K$ and $c_1',\dots, c_K'$, where each $c_i'$ is parallel to $c_i$ in $S_n$.  As before, we may assume each $c_i$ is an essential curve in $S_n$.  Let $d_i$ (resp. $d_i'$) be the disk in $\Delta$ bounded by $c_i$ (resp. $c_i'$).  We use $P_\Delta$ to denote the closure of $\Delta-\cup_{i=1}^k\Delta_i-\cup_{i=1}^Kd_i-\cup_{i=1}^Kd_i'$.  So $P_\Delta$ is a properly embedded planar surface in the handlebody $H_1$, and by our previous assumptions, each component of $\partial P_\Delta$ is essential in $S_n$.  By Lemma~\ref{Ls22}, each circle in $\partial P_\Delta$ bounds a compressing disk in a handlebody.  Since each $\partial\Delta_i$ bounds a disk in $H_2$ and the Heegaard surface $S_n$ is strongly irreducible, each component of $\partial P_\Delta$ bounds a compressing disk in $H_2$.  By Corollary~\ref{Csch}, if $P_\Delta$ is incompressible in $H_1$, then $P_\Delta$ is $\partial$--parallel in $H_1$. 

Similar to the arguments for the innermost monogons, we can take $2N$ parallel copies of $E$ and use $N$ disjoint HTH bands to construct $N$ pinched disks, $\tilde{\Delta}_1, \dots, \tilde{\Delta}_N$.  Since these pinched disks are constructed using parallel copies of the same monogon $E$, we may apply the arguments for $\Delta$ and $P_\Delta$ above to each of the $N$ pinched disks $\tilde{\Delta}_1, \dots, \tilde{\Delta}_N$.  Let $P_1,\dots,P_N$ be the planar sub-surfaces of these $N$ pinched disks constructed in the same way as the $P_\Delta$ above.  In particular, each $P_i$ is properly embedded in $H_1$ and each circle in $\partial P_i$ bounds a compressing disk in $H_2$.  Each boundary circle of $P_i$ is either the boundary of a pinched disk or a circle parallel to some $c_j$ in $S_n$.  To simplify notation, we assume each $P_i$ is incompressible.  The proof for the compressible case is the same after we compress the $P_i$'s into incompressible pieces, as in Subcase 1b.  So, by Corollary~\ref{Csch}, each $P_i$ is $\partial$--parallel in $H_1$.

Let $W$ be the 3--manifold obtained by adding $K$ 2--handles to $H_1$ along these $c_1,\dots, c_K$.  Since the $N$ pinched disks are constructed using parallel copies of the same monogon $E$, each $P_i$ can be extended to a properly embedded planar surface $\hat{P}_i$ in $W$.  $\hat{P}_i$ can be considered as the planar surface obtained by capping off the $c_i$'s and $c_i'$'s by disks.  So, by our assumption on $\partial P_i$, each boundary circle of $\hat{P}_i$ is the boundary of a pinched disk which is either some $\tilde{\Delta}_j$ or a pinched disk in $int(\tilde{\Delta}_j)$ formed by innermost monogons.  

By the construction in section~\ref{Shelix}, there is a rectangle in $S_n$ with two opposite edges glued to the boundary of each pinched disk, as shown in the shaded regions in Figure~\ref{monogon} (b).  Since the $c_i$'s are circles in $E$, these rectangles are in $\partial W$.  Hence there is such a rectangle in $\partial W$ attached to each boundary circle of $\hat{P}_i$.  Moreover, for any two adjacent pinched disks, there is also a short $\eta$--arc connecting the two rectangles, as shown in Figure~\ref{monogon} (b).  Similar to the argument on innermost monogons, we may assume these $\tilde{\Delta}_i$'s have a natural order in the following sense:  If $R_i$ is a rectangle attached to a circle in $\partial\hat{P}_i$ with $2\le i\le N-1$, then as shown in Figure~\ref{monogon} (c), there are two arcs $\eta_{i-1}$ and $\eta_i$ connecting $R_i$ to two rectangles $R_{i-1}$ and $R_{i+1}$, where $R_{i-1}$ (resp. $R_{i+1}$) is a rectangle attached to a circle in $\partial\hat{P}_{i-1}$ (resp. $\partial\hat{P}_{i+1}$).  Therefore, we may assume that, if $i\ne 1$ and $i\ne N$, there are two $\eta$--arcs for each component of $\partial\hat{P}_i$, connecting the attached rectangle to $\partial\hat{P}_{i-1}$ and $\partial\hat{P}_{i+1}$, as shown in Figure~\ref{monogon} (c), where $\hat{P}_{i-1}$ and $\hat{P}_{i+1}$ are different planar surfaces.  The fact that $\hat{P}_{i-1}$ and $\hat{P}_{i+1}$ are different surfaces is important to our proof.

Since each $P_i$ is $\partial$--parallel in $H_1$, each $\hat{P}_i$ must be $\partial$--parallel in $W$.  Let $Q_i\subset\partial W$ be the sub-surface of $\partial W$ that is parallel to $\hat{P}_{i}$ and with $\partial Q_i=\partial\hat{P}_{i}$.  Since these $\hat{P}_i$'s are disjoint, any two planar surfaces $Q_i$ and $Q_j$ are either disjoint or nested.

Similar to the arguments on the innermost monogons, let $m_1,\dots,m_{2K}$ be the $2K$ disks in $\partial W-S_n$ parallel to the cores of the 2--handles added to $H_1$. We first suppose some $Q_k$ ($1\le k\le N$) does not contain any disk $m_i$.  Since any planar surface inside $Q_k$ does not contain any disk $m_i$ either, we may assume $Q_k$ is innermost,   Then by the arguments in Subcase 1b on the $Q_i$'s, there must be a rectangle $R$ attached to $\partial Q_k$ and lying in $S_n-int(Q_k)$.  So the $\eta$--arc attached to $R$ must lie in $Q_k$ and hence $Q_k$ must contain another planar surface $Q_j$ ($j\ne k$), which contradicts the assumption that $Q_k$ is innermost.  Thus, we may assume each $Q_k$ contains some disk $m_i$.

Since $K$ is bounded by a number independent of $S_n$, similar to the arguments on the innermost monogons above, if $N$ is large enough, we can find 3 nested planar surfaces, say $Q_{n_1}\subset Q_{n_2}\subset Q_{n_3}$, such that $Q_{n_3}-Q_{n_2}$ and $Q_{n_2}-Q_{n_1}$ do not contain any disk $m_i$.  Moreover, if $N$ is large, we can find many such nested planar surfaces so that $n_2\ne 1$ and $n_2\ne N$.  Since each $Q_k$ contains some disk $m_i$, $Q_{n_3}-Q_{n_2}$ and $Q_{n_2}-Q_{n_1}$ do not contain any other planar surface $Q_k$.  Moreover, we can choose the $Q_{n_1}$, $Q_{n_2}$ and $Q_{n_3}$ so that there is no $Q_k$ satisfying $Q_{n_1}\subset Q_{k}\subset Q_{n_2}$ or $Q_{n_2}\subset Q_{k}\subset Q_{n_3}$.

Let $\alpha$ be a boundary circle of $Q_{n_2}$.  Since $Q_{n_3}$ is a planar surface and $Q_{n_2}\subset Q_{n_3}$, $\alpha$ is separating in $Q_{n_3}$ and bounds a sub-surface $Q_\alpha$ in $Q_{n_3}$.  We can choose $\alpha$ so that $Q_{n_2}\subset Q_\alpha$. Let $R$ be the rectangle attached to this boundary circle $\alpha$ of $Q_{n_2}$.  By our assumption on $n_2$, there is a pair of $\eta$--arcs connecting the rectangle $R$ to two different planar surfaces.  However, by our assumptions on $Q_{n_1}$, $Q_{n_2}$, $Q_{n_3}$ and $\alpha$, if $R\subset Q_{n_2}$, both $\eta$--arcs must connect $R$ to $\partial Q_{n_3}$; if $R\subset S_n-int(Q_{n_2})$, both $\eta$--arcs must connect $R$ to $\partial Q_{n_1}$, which contradicts previous assumption that the pair of $\eta$--arcs connect $R$ to different $\hat{P}_i$'s, see Figure~\ref{monogon} (c).  This finishes the proof Lemma~\ref{Lincomp}.
\end{proof}

\begin{lemma}\label{Lendincomp}
$\mu$ is end-incompressible.
\end{lemma}
\begin{proof}[Proof of Lemma~\ref{Lendincomp}]
As before, by Proposition~\ref{Psplit}, we can split $B^-$ and $B$ so that $B^-$ has no disk of contact and fully carries $\mu$.  We may also split $B^-$ so that the number of components of $M-B^-$ is the smallest among all the branched surfaces fully carrying $\mu$.  After some isotopy, we may assume that $\partial_hN(B^-)\subset\mu$.  Since $\mu$ is incompressible by Lemma~\ref{Lincomp}, $\partial_hN(B^-)$ is incompressible in $M-int(N(B^-))$.  Suppose $\mu$ is not end-incompressible and let $E$ be a monogon in $M-int(N(B^-))$.  Let $\hat{E}$ be the component of $M-int(N(B^-))$ containing $E$.  By Proposition~\ref{Pmon}, $\hat{E}$ must be a solid torus of the form $E\times S^1$.  Let $L$ be the leaf that contains the horizontal boundary component of $\hat{E}$.  Since $|M-B^-|$ is the smallest, we cannot split $N(B^-)$ along $L$ connecting $\hat{E}$ to other components of $M-int(N(B^-))$.   

We may assume $L$ is an orientable surface.  We claim that $L$ must be an infinite annulus.  If $L$ is not an infinite annulus, we can construct a compressing disk for $L$ by connecting two parallel copies of the monogon $E$ and a long vertical band, as shown in Figure~\ref{monogon} (a), similar to the construction of a pinched disk before.  Thus $L$ is an infinite annulus.  Since $B^-$ does not carry any 2--sphere or torus, this contradicts Lemma~\ref{Lleaf}.
\end{proof}

Since $B^-$ does not carry any 2--sphere, Lemmas~\ref{Lincomp} and \ref{Lendincomp} imply that $\mu$ is an essential lamination.  By Proposition~\ref{PHaken}, $M$ is Haken, which contradicts the hypothesis.  This finishes the proof for part A.

\medskip

\noindent\underline{\textbf{Part B}}.  $\mu$ consists of compact leaves.

\medskip

The only difference between the proofs for Part A and Part B is the construction of $P\times I$.  By Theorem~\ref{TMS}, we may assume $\mu$ is either a family of parallel orientable closed surfaces or a twisted family of parallel closed surfaces.  In both cases, $\mu$ corresponds to a rational point in $\mathcal{PL}(B)$.  For any non-orientable surface $S$ carried by $B$, the boundary of a twisted $I$--bundle over $S$ is an orientable closed surface carried by $B$ and corresponding to the same point in $\mathcal{PL}(B)$ as $S$.  Thus, by using the boundary of a twisted $I$--bundle if necessary, we may assume $\mu$ consists of orientable closed surfaces.  Let $B^-$ be the sub-branched surface of $B$ fully carrying $\mu$.  By Proposition~\ref{Psplit}, after some splittings, we may assume $B^-$ is an orientable closed surface and $N(B^-)$ is a product of an interval and the closed surface $B^-$. Moreover, by Corollary~\ref{Cnormal} and our assumptions on $B$ before, $B^-$ is a normal surface in $M$ with genus at least 2.

We first prove that there must be a non-separating simple closed curve in $B^-$ that bounds an embedded disk $D$ in $M$ (note that $int(D)\cap B^-$ may not be empty).  Since $M$ is non-Haken, $B^-$ is compressible and we can perform a compression on $B^-$ and get a new surface which must also be compressible.  So we can successively perform compressions on the resulting surfaces until we get a collection of $2$--spheres.  If the boundary circle of every compressing disk is separating, then after some compressions, we get an embedded torus.  As every essential simple closed curve in a torus is non-separating, we get a non-separating simple closed curve $\gamma$ in $B^-$ such that $\gamma$ bounds an embedded disk $D$ in $M$.  Moreover, we may assume that $D$ is transverse to $B^-$ and every component of $int(D)\cap B^-$ is a separating curve in $B^-$. 

Let $\gamma_1$ and $\gamma_2$ be two parallel copies of $\gamma$ in $B^-$.  Each $\gamma_i$ bounds a disk $D_i$ in $M$ ($i=1,2$), and each $D_i$ is parallel to $D$. We may assume $D_1\cap D_2=\emptyset$.  Since $\gamma$ is non-separating, there is an arc $\alpha\subset B^-$ connecting $\gamma_1$ to $\gamma_2$, forming a graph $\Gamma=\gamma_1\cup\alpha\cup\gamma_2$, such that $B^--\Gamma$ contains no disk component.  Moreover, since every component of $int(D)\cap B^-$ is a separating curve in $B^-$, we may choose $\alpha$ so that $\alpha\cap int(D_i)=\emptyset$.   Let $A_i$ ($i=1,2$) be an annular neighborhood of $\gamma_i$ in $B^-$ and $Q$ a small neighborhood of $\alpha$ in $B^-$.  Then $P=A_1\cup Q\cup A_2$ is a sub-surface of $B^-$ and no boundary circle of $P$ bounds a disk in $B^-$.  Let $P\times I=\pi^{-1}(P)$ and $A_i\times I=\pi^{-1}(A_i)$ ($i=1,2$), where $\pi: N(B^-)\to B^-$ is the projection.  We may consider $P$ as the limit of $\{S_n\cap(P\times I)\}$ in the corresponding projective lamination space.  We will use this $P\times I$ to construct our HTH bands, as in section~\ref{Shelix}.  

As before, we may assume the sequence of surfaces $\{S_n\}$ satisfy the hypotheses of Lemma~\ref{Lcor}.  By Lemma~\ref{Lessential}, we may assume $S_n\cap (\gamma_i\times I)$ ($i=1,2)$ does not contain any circle that is trivial in $S_n$, for each $n$.  If $S_n\cap(\gamma_i\times I)$ consists of circles, then each circle is essential in $S_n$ and hence bounds a compressing disk in one of the two handlebodies by Lemma~\ref{Ls22}.  However, if $S_n\cap(\gamma_i\times I)$ consists of circles, by Corollary~\ref{Clinear} and Example~\ref{exspiral},  the number of circles in $S_n\cap(\gamma_i\times I)$ tends to infinity as $n$ goes to infinity.  This gives a contradiction to Lemma~\ref{Lbound}.  So $S_n\cap(\gamma_i\times I)$ cannot be a union of circles if $n$ is large enough.  Hence we may assume $S_n\cap(\gamma_i\times I)$ consists of spirals for each $n$.

Therefore, after splitting $B$, we may assume $S_n\cap(A_i\times I)$ ($i=1,2)$ consists of spiraling disks and $S_n\cap (P\times I)$ satisfies the conditions in Example~\ref{exjoint}.  We use the same notations as section~\ref{Shelix}.  Let $h_i$ be the number of components of $S_n\cap(A_i\times I)$ ($i=1,2$), and we may assume $n$ is sufficiently large.  Since $\gamma_1$ and $\gamma_2$ are parallel in $B^-$, we may assume $h_1=h_2$.  Then we can use Example~\ref{HTH} to construction our HTH bands and the remainder of the proof is the same as Part A.  This finishes the proof of Theorem~\ref{Tfinite} and Theorem~\ref{main}.
\end{proof}

\section{The Casson-Gordon example}\label{Sexample}

Casson and Gordon gave an example of a 3--manifold that has an infinite family of strongly irreducible Heegaard splittings with different genera \cite{CG2}, see \cite{Ko, Sed}.  By Theorem~\ref{main}, such a 3--manifold must be Haken.  In  fact, it is easy to directly show that the 3--manifolds in the Casson-Gordon example are Haken.  The proof of Theorem~\ref{main} indicates that there should be an incompressible surface as the limit of the infinite family of Heegaard surfaces.  In this section, we construct such an incompressible surface.

\begin{figure}
\begin{center}
\includegraphics{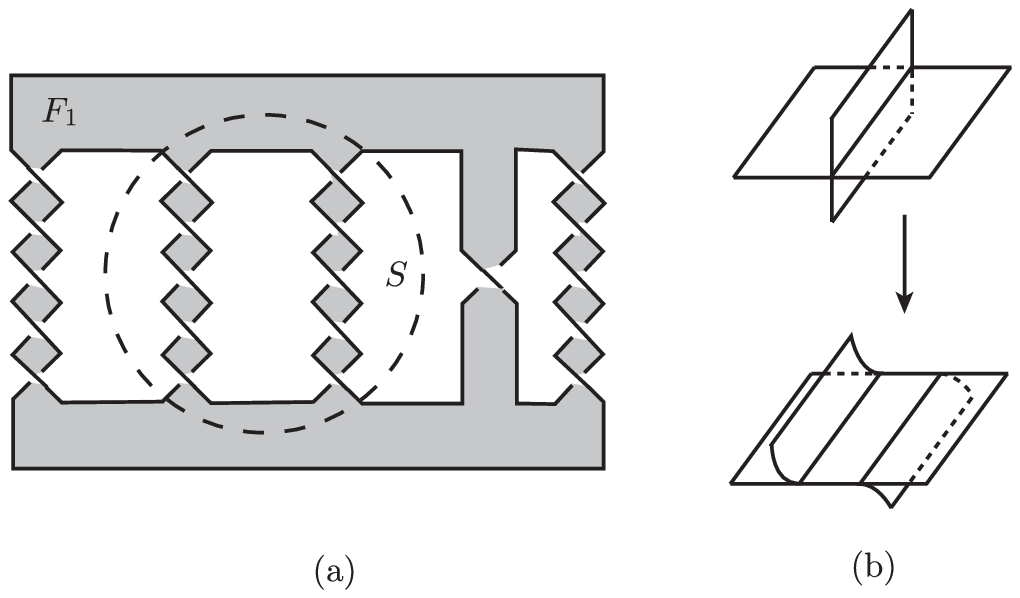}
\caption{}\label{pretzel}
\end{center}
\end{figure}

Before carrying out the construction, we give a brief overview of the Casson-Gordon example and we refer to \cite{Ko, Sed} for more details.  We first take a pretzel knot $K=(p_1,p_2,p_3,1,p_4)$ in $S^3$, where $|p_i|\ge 5$.  The standard Seifert surface $F_1$ from the Seifert algorithm is a free Seifert surface.  Let $S$ be a 2--sphere in $S^3$ that cuts the knot into 2 tangles, as shown in Figure~\ref{pretzel} (a).  If we flip a tangle bounded by $S$ along a horizontal axis by $180^\circ$, we get the same knot with a different projection $(p_1,-1,p_2,p_3,1,1,p_4)$. By a theorem of Parris \cite{P}, the standard Seifert surface $F_2$ from the Seifert algorithm is also a free Seifert surface with $genus(F_2)=genus(F_1)+1$.  By flipping the tangle $k$ times, we get an infinite family of free Seifert surfaces $\{F_k\}$ with increasing genus.  

Let $\eta(K)$ be a tubular neighborhood of the knot $K$ and let $M_0=S^3-\eta(K)$ be the knot exterior.  Let $H_k$ be the closure of a small neighborhood of $F_k$ in $M_0$.  So $H_k$ is a handlebody.  Since $F_k$ is a free Seifert surface, $\overline{M_0-H_k}$ is also a handlebody.  Let $K(p/q)$ be the closed manifold obtained by the Dehn filling to $M_0$ along the slope $p/q$.  We may regard $H_k$ as a handlebody in $K(p/q)$.  In fact, if $p=1$, $K(1/q)-int(H_k)$ is also a handlebody and $S_k=\partial H_k$ is a Heegaard surface for $K(1/q)$.  Casson and Gordon showed that \cite{CG2, Ko, MSch}, if $|q|\ge 6$, then this Heegaard splitting of $K(1/q)$ by $S_k=\partial H_k$ is strongly irreducible.  So we get an infinite family of strongly irreducible Heegaard surfaces $\{S_k\}$ for $M=K(1/q)$ ($|q|\ge 6$).

In \cite{Ko}, Kobayashi gave an interpretation of the sequence of free Seifert surfaces $\{F_k\}$ through branched surfaces.  Let $F_1$ be the free Seifert surface of $M_0=S^3-\eta(K)$ above, and $S$ the punctured 2--sphere as shown in Figure~\ref{pretzel}~(a).  By fixing a normal direction for $F_1$ and $S$, we can deform $F_1\cup S$ into a branched surface $B_0$, as shown in Figure~\ref{pretzel} (b).  Both $F_1$ and $S$ are carried by $B_0$, so we can assume $F_1$ and $S$ lie in  $N(B_0)$, a fibered neighborhood of $B_0$.  Then the canonical cutting and pasting on $F_1$ and $S$ produce another Seifert surface $F_1+S$.  Kobayashi showed that $F_2=F_1+S$ is the same free Seifert surface described above.  Moreover, $F_k=F_1+(k-1)S$.

\begin{figure}
\begin{center}
\includegraphics{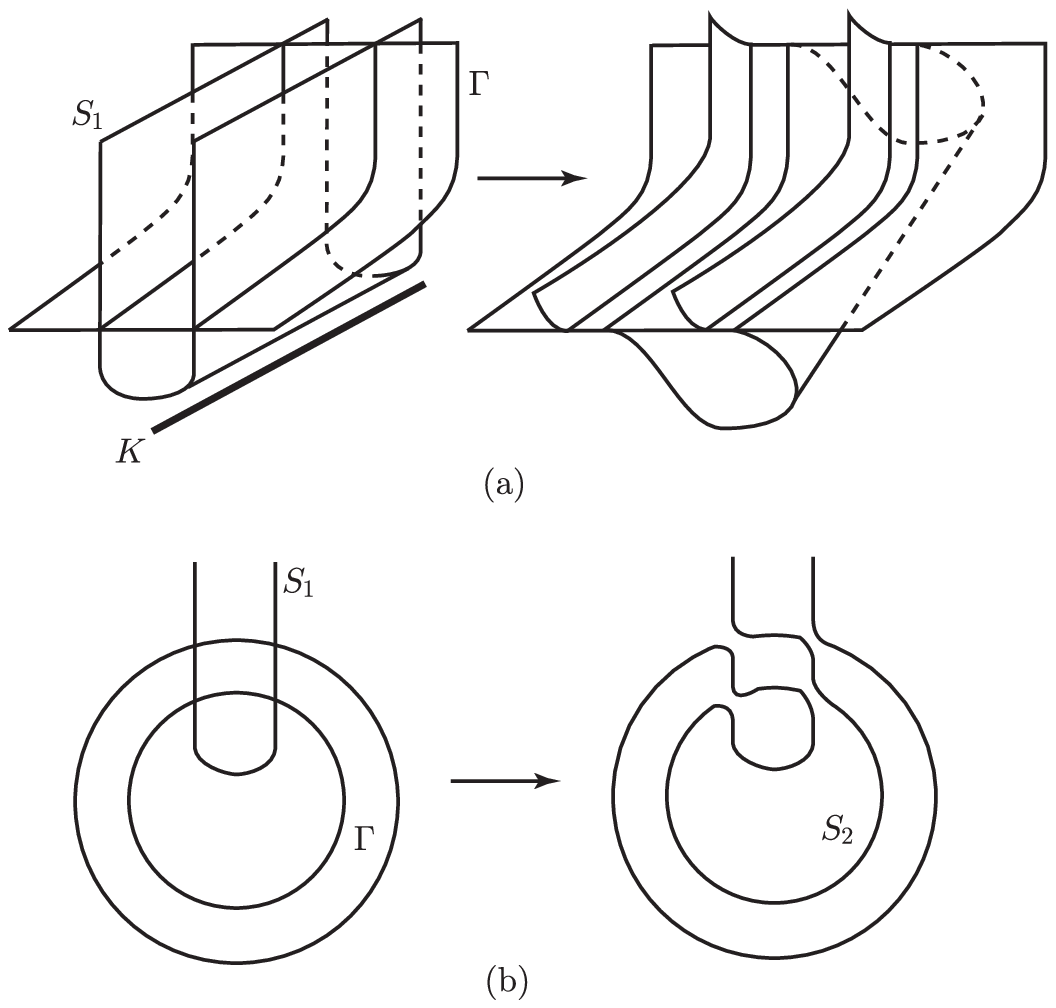}
\caption{}\label{deform}
\end{center}
\end{figure}

As we mentioned before, the closed manifold $M=K(1/q)$ is Haken.  The 2--sphere $S\subset S^3$ in Figure~\ref{pretzel} (a) cuts $(S^3, K)$ into a pair of non-trivial tangles $(E_1,K_1)$ and $(E_2, K_2)$, where $E_1$ and $E_2$ are the pair of 3--balls in $S^3$ bounded by $S$ and $K_i\subset E_i$ is a pair of strings.  Let $\eta(K_1)$ be a small neighborhood of $K_1$ in $E_1$.  Then $\Gamma=\partial(E_1-\eta(K_1))$ is a closed surface of genus 2 in $S^3-K$.  It is not hard to see that $\Gamma$ is incompressible in $S^3-K$ (for instance see \cite{Wu}).  By a theorem of Menasco \cite{Me}, $\Gamma$ remains incompressible after any non-trivial Dehn surgery on $K$.

Next we will show that $\Gamma$ can be considered as the limit of the sequence of Heegaard surfaces $\{S_k\}$.  

We start with the Seifert surface $F_1$ and consider $K=\partial F_1$ ($F_1\cap\Gamma\ne\emptyset$).  Let $\eta(F_1)$ be a small neighborhood of $F_1$ in $S^3$.  After moving $K$ slightly off $\eta(F_1)$, we can regard the Heegaard surface $S_1$ of $M=K(1/q)$ as the boundary surface of the closure of $\eta(F_1)$.  $S_1\cap\Gamma$ consists of closed curves. 

Similar to the construction of the branched surface $B_0$ above, we can deform $S_1\cup\Gamma$ into a branched surface $B$ as shown in Figure~\ref{deform} (a).  $B$ carries both $S_1$ and $\Gamma$, so we can assume $\Gamma$ and $S_1$ lie in $N(B)$ and transverse to the $I$--fibers.  Then we perform the canonical cutting and pasting on $S_1$ and 2 parallel copies of $\Gamma$, as shown in Figure~\ref{deform} (b). It is not hard to see that the resulting surface $S_1+2\Gamma$ is isotopic to $S_2$.  Similarly, $S_3=S_2+2\Gamma$ and $S_k=S_1+2(k-1)\Gamma$.  By our discussion on projective lamination spaces, $\Gamma$ is indeed the limit of the sequence of Heegaard surfaces $\{S_k\}$.


\begin{thebibliography}{50}


\bibitem{AL}
Ian Agol and Tao Li \emph{An algorithm to detect laminar 3--manifolds}. Geom. Topol. 7 (2003) 287--309

\bibitem{BCZ}
M. Boileau and D.J. Collins and H. Zieschang,  \emph{Genus $2$ Heegaard decompositions of small Seifert manifolds}. Ann. Inst. Fourier (Grenoble) \textbf{41} (1991), no. 4, 1005--1024.

\bibitem{BO}
F. Bonahon and J.P. Otal, \emph{Scindements de Heegaard des espaces lenticulaires}, Ann. Sci. Ec. Norm. Sup. 16 (4) (1983) 451--466

\bibitem{CG}
Andrew Casson and Cameron Gordon, \emph{Reducing Heegaard splittings}.  Topology and its Applications, \textbf{27} 275--283 (1987).

\bibitem{CG2}
Andrew Casson and Cameron Gordon, unpublished.

\bibitem{FO} 
William Floyd and Ulrich Oertel,  \emph{Incompressible surfaces via branched surfaces}. Topology \textbf{23} (1984), no. 1, 117--125. 

\bibitem{G7}
David Gabai,  \emph{Foliations and 3--manifolds}. Proceedings of the International Congress of Mathematicians, Vol. I (Kyoto, 1990) 609--619. 

\bibitem{GO}
David Gabai and Ulrich Oertel, 
\emph{Essential laminations in $3$--manifolds}.  Ann. of Math. (2), \textbf{130} (1989) 41--73

\bibitem{H}
Wolfgang Haken, \emph{Some results on surfaces in 3--manifolds}.  Studies in Modern Topology, Math. Assoc. Amer., distributed by Pretice-Hall, (1968) 34--98.

\bibitem{Hat}
Allen Hatcher, \emph{Measured lamination spaces for surfaces, from the topological viewpoint}.  Topology and its Applications \textbf{30} (1988) 63--88.

\bibitem{Jo1}
Klaus Johannson,  \emph{Heegaard surfaces in Haken $3$--manifolds}. Bull. Amer. Math. Soc. \textbf{23} (1990), no. 1, 91--98.

\bibitem{Jo2}
Klaus Johannson,  \emph{Topology and combinatorics of 3--manifolds}. Lecture Notes in Mathematics, \textbf{1599}, Springer-Verlag, Berlin, 1995.

\bibitem{Ko}
Tsuyoshi Kobayashi, \emph{A construction of 3--manifolds whose homeomorphism classes of Heegaard splittings have polynomial growth}.  Osaka Journal of Mathematics, \textbf{29} (1992), 653--674.

\bibitem{La}
Marc Lackenby, \emph{The asymptotic behaviour of Heegaard genus}. Preprint. arXiv: math.GT/0210316 

\bibitem{L1}
Tao Li, \emph{Laminar branched surfaces in 3--manifolds.}  Geometry and Topology, \textbf{6} (2002), 153--194. 

\bibitem{L5}
Tao Li, \emph{Boundary curves of surfaces with the 4-plane property}.  Geometry and Topology,  \textbf{6} (2002), 609--647.

\bibitem{L2}
Tao Li, \emph{An algorithm to find vertical tori in small Seifert fiber spaces}. Comment. Math. Helv. \textbf{81} (2006), 727--753.

\bibitem{L4}
Tao Li, \emph{Heegaard surfaces and measured laminations, I: the Waldhausen conjecture}. Invent. Math.  \textbf{167} (2007) 135--177.

\bibitem{MMZ}
Joseph Masters and William Menasco and Xingru Zhang, \emph{Heegaard splittings and virtually Haken Dehn filling}. Preprint, arXiv:math.GT/0210412

\bibitem{Me}
William Menasco, \emph{Closed incompressible surfaces in alternating knot and link complements}. Topology \textbf{23} (1984), no. 1, 37--44.

\bibitem{MS}
John Morgan and Peter Shalen, \emph{Degenerations of hyperbolic structures, II: Measured laminations in 3--manifolds}. Annals of Math. \textbf{127} (1988), 403--456.

\bibitem{M}
Yoav Moriah,  \emph{Heegaard splittings of Seifert fibered spaces}. Invent. Math. \textbf{91} (1988), 465--481.

\bibitem{MSch}
Yoav Moriah and Jennifer Schultens, \emph{Irreducible Heegaard splittings of Seifert fibered spaces are either vertical or horizontal}.  Topology \textbf{37} (1998) 1089--1112.

\bibitem{MSS}
Yoav Moriah and Saul Schleimer and Eric Sedgwick, \emph{Heegaard splittings of the form H + nK}, to appear in Communications in Analysis and Geometry, arXiv:math.GT/0408002


\bibitem{O}
Ulrich Oertel, \emph{Measured laminations in 3--manifolds}.  Trans. Amer. Math. Soc.  \textbf{305}  (1988),  no. 2, 531--573.

\bibitem{P}
Parris, \emph{Pretzel Knots}.  Ph.D. Thesis, Princeton University (1978).

\bibitem{R}
Hyam Rubinstein, \emph{Polyhedral minimal surfaces, Heegaard splittings and decision problems for 3-dimensional manifolds}. Proc. Georgia Topology Conference, Amer. Math. Coc./Intl. Press, 1993.

\bibitem{RS}
Hyam Rubinstein; Martin Scharlemann, \emph{Comparing Heegaard splittings of non-Haken $3$--manifolds}.  Topology  \textbf{35} (1996),  no. 4, 1005--1026

\bibitem{S}
Martin Scharlemann, \emph{Local detection of strongly irreducible Heegaard splittings}. Topology and its Applications, \textbf{90} (1998) 135--147.

\bibitem{Sed}
Eric Sedgwick, \emph{An infinite collection of Heegaard splittings that are equivalent after one stabilization}.  Math. Ann. \textbf{308}, 65--72 (1997).

\bibitem{St}
Michelle Stocking, \emph{Almost normal surfaces in 3--manifolds}. Trans. Amer. Math. Soc. \textbf{352}, 171--207 (2000).


\bibitem{W1}
Friedhelm Waldhausen, \emph{Heegaard-Zerlegungen der $3$--Sph\"{a}re}. Topology \textbf{7} 1968 195--203.

\bibitem{W2}
Friedhelm Waldhausen, \emph{Some problems on 3--manifolds}, Proc. Symp. Pure Math. 32 (1978) 313--322.

\bibitem{Wu}
Ying-Qing Wu, \emph{Dehn surgery on arborescent knots}. J. Differential Geom. \textbf{43} (1996), no. 1, 171--197.

\end{thebibliography}
\end{document}